\documentclass[12pt]{article}
\usepackage{graphicx}
\usepackage{amsmath,amsthm,amssymb,enumerate}%, esint}
\usepackage{euscript,mathrsfs}
\usepackage{color}
\usepackage{dsfont}
\usepackage[left=2cm,right=2cm,top=3.5cm,bottom=3.5cm]{geometry}
\usepackage{color}
\usepackage[framemethod=tikz]{mdframed}
\allowdisplaybreaks

\usepackage{soul}

\catcode`\@=11 \@addtoreset{equation}{section}

\catcode`\@=12

\allowdisplaybreaks

\newtheorem{Theorem}{Theorem}[section]
\newtheorem{Proposition}[Theorem]{Proposition}
\newtheorem{Lemma}[Theorem]{Lemma}
\newtheorem{Corollary}[Theorem]{Corollary}

\theoremstyle{definition}
\newtheorem{Definition}[Theorem]{Definition}

\newtheorem{Remark}[Theorem]{Remark}

\newcommand{\bTheorem}[1]{
\begin{Theorem} \label{T#1} }
\newcommand{\eT}{\end{Theorem}}

\newcommand{\bProposition}[1]{
\begin{Proposition} \label{P#1}}
\newcommand{\eP}{\end{Proposition}}

\newcommand{\bLemma}[1]{
\begin{Lemma} \label{L#1} }
\newcommand{\eL}{\end{Lemma}}

\newcommand{\bCorollary}[1]{
\begin{Corollary} \label{C#1} }
\newcommand{\eC}{\end{Corollary}}

\newcommand{\bRemark}[1]{
\begin{Remark} \label{R#1} }
\newcommand{\eR}{\end{Remark}}

\newcommand{\bDefinition}[1]{
\begin{Definition} \label{D#1} }
\newcommand{\eD}{\end{Definition}}

\newcommand{\Del}{\Delta_x}

\newcommand{\Ds}{\mathbb{D}_x}

\DeclareMathOperator{\supp}{supp}

\newcommand{\tS}{\widetilde{S}}
\newcommand{\bfphi}{\boldsymbol{\varphi}}

\newcommand{\bFormula}[1]{
\begin{equation} \label{#1}}
\newcommand{\eF}{\end{equation}}

\newcommand{\Cd}{c({\rm data})}

\newcommand{\Ov}[1]{\overline{#1}}

\newcommand{\aleq}{\stackrel{<}{\sim}}

\newcommand{\vr}{\varrho}

\newcommand{\tvr}{\tilde \vr}
\newcommand{\tvu}{{\tilde \vu}}
\newcommand{\tvt}{\tilde \vt}

\newcommand{\vt}{\vartheta}
\newcommand{\vu}{\vc{u}}
\newcommand{\vm}{\vc{m}}

\newcommand{\Cdk}{c({\rm data}, K)}
\newcommand{\vc}[1]{{\bf #1}}

\newcommand{\Div}{{\rm div}_x}
\newcommand{\Grad}{\nabla_x}

\newcommand{\dx}{\,{\rm d} {x}}

\newcommand{\dt}{\,{\rm d} t }

\newcommand{\intO}[1]{\int_{\Omega} #1 \ \dx}

\newcommand{\D}{{\rm d}}

\newcommand{\ep}{\varepsilon}

\newcommand{\vtB}{\vt_B}

\newcommand{\br}{ \nonumber \\ }

\newcommand{\ThB}{\Theta_{\rm B}}
\newcommand{\ThU}{\Theta_{\rm U}}
\def\softd{{\leavevmode\setbox1=\hbox{d}%
          \hbox to 1.05\wd1{d\kern-0.4ex{\char039}\hss}}}
\definecolor{Cgrey}{rgb}{0.85,0.85,0.85}
\definecolor{Cblue}{rgb}{0.50,0.85,0.85}
\definecolor{Cred}{rgb}{1,0,0}
\definecolor{fancy}{rgb}{0.10,0.85,0.10}

\newcommand\Cbox[2]{%
    \newbox\contentbox%
    \newbox\bkgdbox%
    \setbox\contentbox\hbox to \hsize{%
        \vtop{
            \kern\columnsep
            \hbox to \hsize{%
                \kern\columnsep%
                \advance\hsize by -2\columnsep%
                \setlength{\textwidth}{\hsize}%
                \vbox{
                    \parskip=\baselineskip
                    \parindent=0bp
                    #2
                }%
                \kern\columnsep%
            }%
            \kern\columnsep%
        }%
    }%
    \setbox\bkgdbox\vbox{
        \color{#1}
        \hrule width  \wd\contentbox %
               height \ht\contentbox %
               depth  \dp\contentbox
        \color{black}
    }%
    \wd\bkgdbox=0bp%
    \vbox{\hbox to \hsize{\box\bkgdbox\box\contentbox}}%
    \vskip\baselineskip%
}

\mdfdefinestyle{MyFrame}{%
	linecolor=black,
	outerlinewidth=1pt,
	roundcorner=5pt,
	innertopmargin=\baselineskip,
	innerbottommargin=\baselineskip,
	innerrightmargin=10pt,
	innerleftmargin=10pt,
	backgroundcolor=white!20!white}

\date{}

\begin{document}

%%%%%%%%%%%%%%%%%%%%%%%%%%%%%%%%

\title{The Rayleigh--B\' enard problem for compressible fluid flows}

\author{Eduard Feireisl
\thanks{The work of E.F. was partially supported by the
Czech Sciences Foundation (GA\v CR), Grant Agreement
18--05974S. The Institute of Mathematics of the Academy of Sciences of
the Czech Republic is supported by RVO:67985840. This work is partially supported by the Simons Foundation Award No 663281 granted to the Institute of Mathematics of the Polish Academy of Sciences for the years 2021-2023.} \and Agnieszka \'Swierczewska-Gwiazda
\thanks{The work of A. \'S-G. was partially supported by  National Science Centre
(Poland),  agreement no 2017/27/B/ST1/01569.}
}

\date{}

\maketitle

\medskip

\centerline{Institute of Mathematics of the Academy of Sciences of the Czech Republic}

\centerline{\v Zitn\' a 25, CZ-115 67 Praha 1, Czech Republic}

\medskip

\centerline{Institute of Applied Mathematics and Mechanics, University of Warsaw}

\centerline{Banacha 2, 02-097 Warsaw, Poland}

\begin{abstract}
	
	We consider the physically relevant fully compressible setting of the Rayleigh--B\' enard problem of 
	a fluid confined between two parallel plates, heated from the bottom, and subjected to the gravitational force. 
	Under suitable restrictions imposed on the constitutive relations we show that this open system is dissipative in the sense of Levinson, meaning there exists a bounded absorbing set for any global--in--time weak solution. In addition, global--in--time trajectories are asymptotically compact in suitable topologies and the system possesses a global compact 
	trajectory attractor $\mathcal{A}$. The standard technique of Krylov and Bogolyubov then yields the existence of an invariant measure -- a stationary statistical solution sitting on $\mathcal{A}$. In addition, the Birkhoff--Khinchin ergodic theorem 
	provides convergence of ergodic averages of solutions belonging to $\mathcal{A}$ a.s. with respect to the invariant measure.

\end{abstract}

%\bigskip


%{\bf MSC:}
%\bigskip

\tableofcontents

\section{Introduction}
\label{i}

The Rayleigh--B\' enard problem concerns the motion of a fluid confined between two parallel planes, where the temperature of the
bottom plane is maintained at the level $\ThB$, while the top plane has the ambient temperature $\ThU$, typically $\ThU \leq \ThB$. The only volume force is due to gravitation acting in the downward vertical direction. The fluid mass density $\vr = \vr(t,x)$, the temperature 
$\vt(t,x)$, and the velocity $\vu = \vu(t,x)$ obey the standard system of equations of continuum fluid mechanics:

\begin{mdframed}[style=MyFrame]
\begin{align}
	\partial_t \vr + \Div (\vr \vu) &= 0, \label{i1} \\
	\partial_t (\vr \vu) + \Div (\vr \vu \otimes \vu) + \Grad p (\vr, \vt) &= \Div \mathbb{S} + \vr \Grad G, \label{i2}\\
	\partial_t (\vr e(\vr, \vt)) + \Div (\vr e(\vr, \vt) \vu) + \Grad \vc{q} &= \mathbb{S} : \Ds \vu - p(\vr, \vt) \Div \vu,
	\label{i3}
	\end{align}
\end{mdframed}

\noindent where we have denoted

the pressure \dotfill $p = p(\vr,\vt)$,

the (specific) internal energy \dotfill $e = e(\vr, \vt)$,

the viscous stress tensor \dotfill $\mathbb{S}$, 

the gravitational potential \dotfill $G = - gx_3$,

the heat flux \dotfill $\vc{q}$

the symmetric part of the velocity gradient \dotfill $\Ds \vu = \frac{1}{2} \left( \Grad \vu + \Grad^t \vu \right)$.

\medskip

For the sake of simplicity, we consider the periodic boundary conditions with respect to the horizontal variables. Accordingly, the fluid domain $\Omega$ can be identified with 
\begin{equation} \label{i4}
	\Omega = \mathbb{T}^2 \times (0,1),
	\end{equation}
where $\mathbb{T}^2$ is the two--dimensional flat torus. If the boundary planes are impermeable and the viscous fluid sticks to them, the relevant boundary conditions read 
\begin{mdframed}[style=MyFrame]
	
	\begin{align}
	\vu|_{x_3 = 0} &= \vu|_{x_3 = 1} = 0, \label{i5}\\ 
	\vt|_{x_3 = 0} &= \ThB,\ \vt|_{x_3 = 1} = \ThU. \label{i6}
\end{align}	
\end{mdframed} 	

We suppose that the fluid is Newtonian, with the viscous stress 
\begin{equation} \label{i7}
	\mathbb{S} (\vt, \Ds \vu) = \mu(\vt) \left( \Grad \vu + \Grad^t \vu - \frac{2}{3} \Div \vu \mathbb{I} \right) + 
	\eta(\vt) \Div \vu \mathbb{I}.
	\end{equation}
The heat flux is given by Fourier's law 
\begin{equation} \label{i8}
	\vc{q}(\vt, \Grad \vt) = - \kappa(\vt) \Grad \vt,
	\end{equation}
	where  $\kappa$ is the conductivity.
The field equations \eqref{i1}--\eqref{i3}, endowed with the constitutive relations \eqref{i7}, \eqref{i8}, will be referred to 
as \emph{Navier--Stokes--Fourier system}.

The behaviour of the fluid under the boundary conditions \eqref{i5}, \eqref{i6} with a sufficiently large difference between the 
top and bottom temperatures is a prominent example of turbulence, see e.g. Davidson \cite{DAVI}. There is a large amount of mathematical literature devoted to the asymptotic behaviour of solutions to the Rayleigh--B\' enard problem in the simplified incompressible framework, where the system \eqref{i1}--\eqref{i3} is replaced by the Oberbeck--Boussinesq approximation, see Constantin et al. 
\cite{CoFoMaTe}, Foias, Manley and Temam \cite{FoMaTe}, Cao et al. \cite{CaJoTuWh}
and the references therein. Recently, the problem motivated a series of studies by Otto et al. \cite{ChoNobOtt},
\cite{NobOtt1}, \cite{OTTOetal1}
concerning 
the associated scaling laws.

Much less seems to be known in the original and physically relevant framework of compressible and heat conducting fluids. Here, the rigorous analysis is hampered by the following notoriously known difficulties: 

\begin{itemize}
	\item Navier--Stokes--Fourier system endowed with the boundary conditions \eqref{i5}, \eqref{i6} is an open system in the regime far from equilibrium. Unfortunately, the existence of \emph{global--in--time} smooth solutions is known only in the case 
	of closed systems approaching the equilibrium solution in the long run, see Matsumura and Nishida 
	\cite{MANI1}, \cite{MANI}, or Valli and Zajaczkowski 
	\cite{VAZA}.
	
	\item The available theory of weak solutions developed in \cite{FeNo6A} (see also the alternative approach by 
	Bresch and Desjardins \cite{BRDE} and Bresch and Jabin \cite{BreJab}) applies to conservative or periodic boundary conditions pertinent to the closed systems. Note that the dynamics of the Navier--Stokes--Fourier system with 
	\emph{conservative} boundary conditions is nowadays well understood, see \cite{FP20}. Indeed, in accordance with the celebrated statement of Clausius:
	
	{\it ``Die Energie der Welt ist konstant. Die Entropie der Welt strebt einem Maximum zu'' }

	\rightline{{\it Rudolf Clausius}, Poggendorff's Annals of Physics 1865 (125), 400;}
	
	any global--in--time weak solution of the Navier--Stokes--Fourier system with conservative boundary conditions and driven by a conservative volume force tends to an equilibrium, see e.g. \cite{FP9}, \cite{FP20} or Novotn\' y and Pokorn\' y \cite{NovPok07}, Novotn\' y and Stra\v skraba \cite{NOS1}, \cite{NOST}. 
	Here, conservative means that the system is thermally insulated, the boundary conditions \eqref{i6} being replaced by 
	\[
	\vc{q} \cdot \vc{n} = q_3 = 0 \ \mbox{on}\ \partial \Omega.
	\] 
	
	\item The weak solutions are not  (known to be) uniquely determined by the initial/boundary data.

	\end{itemize}

Recently, the theory of weak solutions has been extended to non--zero in/out flux boundary conditions in \cite{FeiNov20}, and, finally, to general Dirichlet boundary conditions in \cite{ChauFei}. In particular, the theory of weak solutions proposed in 
\cite{ChauFei} yields a suitable platform to attack the Rayleigh--B\' enard problem \eqref{i5}, \eqref{i6}. As pointed out above, the weak solutions are not known to be uniquely determined by the initial/boundary data. Accordingly, we follow the approach advocated by Sell \cite{SEL} and M\' alek and Ne\v cas \cite{MANE} replacing the standard phase space by the space of trajectories.

The principal objective of the paper is to establish the following two basic results:

\begin{itemize}
	
	\item {\bf Levinson dissipativity or bounded absorbing set.} Any global--in--time weak solution to the Navier--Stokes--Fourier system endowed with the boundary conditions \eqref{i5}, \eqref{i6} enters eventually a bounded absorbing set. In comparison with \cite{FeiKwo}, we relax the hypothesis of the hard sphere pressure and consider the physically relevant equation of 
	state of a general monoatomic gas with the effect of radiation proposed in \cite{FeNo6A}.
	
	\item {\bf Asymptotic compactness.} Similarly to \cite[Chapter 4, Theorem 4.2]{FeiPr}, we show that any bounded family of global solutions is precompact in a suitable topology of the trajectory space, whereas any of its accumulation points represents a weak solution of the same problem.

	\end{itemize}

Using the above results, we establish the existence of a compact trajectory attractor, an invariant 
measure and the existence of stationary statistical solutions generated by bounded trajectories.
Finally, we also discuss the existence of the ergodic averages in the spirit of \cite{FanFeiHof}. 

As pointed out, the key point of the analysis 
is the Levinson dissipativity or the existence of a universal bounded absorbing set for the ``monoatomic'' equation of 
state introduced in \cite[Chapters 1,2]{FeNo6A}. This is rather surprising as this constitutive equation can be seen as a 
temperature dependent counterpart of the isentropic pressure law $p(\vr) \approx \vr^{\gamma}$, with $\gamma = \frac{5}{3}$. Note that 
for the isentropic model, the existence of a bounded absorbing set is known only if $\gamma > \frac{5}{3}$, see \cite{FP14},  
whereas the limit case $\gamma = \frac{5}{3}$ requires smallness of the total mass of the fluid, see Wang and Wang \cite{WanWan}. Moreover, uniform boundedness of global trajectories for the Navier--Stokes--Fourier system 
is a delicate issue. As is known, see \cite{FP20}, \cite{FeiPr}, the energy of \emph{all} global--in--time solutions tends to infinity with growing time as soon as the system is energetically closed and driven by a non--conservative volume force.

Similarly to the incompressible Navier--Stokes system with conservative boundary conditions 
studied by M\' alek and Ne\v cas \cite{MANE} and Sell \cite{SEL}, the large time asymptotic behaviour of 
solutions to the Rayleigh--B\' enard problem is captured by the set of \emph{entire} trajectories $\mathcal{A}$ defined for all $t \in R$ and with 
uniformly bounded total energy and mass. We show that the set $\mathcal{A}$ is {\bf (i)} non--empty, {\bf (ii)} time 
shift invariant, and {\bf (iii)} compact if endowed by a suitable metrics. The standard Krylov--Bogolyubov technique then yields 
the existence of an invariant measure supported in $\mathcal{A}$ -- a stationary statistical solution of the Rayleigh--B\' enard problem. Moreover, the standard Birkhoff--Khinchin ergodic theorems yields the convergence of the ergodic averages 
a.s. with respect to the invariant measure. Uniqueness of the invariant measure for solutions with the same total mass remains an outstanding open problem.

The paper is organized as follows. In Section \ref{w}, we recall the principal constitutive hypotheses and introduce the concept of weak solution. The main results -- the existence of a bounded absorbing set and asymptotic compactness of global--in--time solutions -- are stated in Section \ref{M}. In Section \ref{B}, we show the existence of a bounded absorbing set in terms of the total energy. The implications of the main results on the long--time behavior of the system are discussed in Section~\ref{L}.

\section{Principal hypotheses, weak solutions}
\label{w}

Following \cite{ChauFei} we introduce the concept of weak solution to the Navier--Stokes--Fourier system based on the combination 
of the balance equations for the entropy and the \emph{ballistic energy}
\[
E_{\tvt} (\vr, \vt, \vu) = E(\vr, \vt, \vu) - \tvt \vr s(\vr, \vt),\ 
E(\vr, \vt, \vu) = \frac{1}{2} \vr |\vu|^2 + \vr e(\vr, \vt),
\]
where $s$ is the entropy related to the other thermodynamic functions through Gibbs' equation 
\begin{equation} \label{w1}
	\vt D s(\vr, \vt) = D e(\vr, \vt) + p D \left( \frac{1}{\vr} \right),
	\end{equation}
and $\tvt$ is an arbitrary continuously differentiable function of $(t,x)$ satisfying 
\begin{equation} \label{w2}
	\tvt > 0,\ \tvt|_{x_3 = 0} = \ThB,\ \tvt|_{x_3 = 1} = \ThU.
	\end{equation}
For the sake of simplicity, we suppose that $\ThB$, $\ThU$ are positive \emph{constants}. A generalization of the results 
of the present paper to the 
space or even time dependent boundary temperatures is possible with obvious modifications in the proofs.

\subsection{Weak solution}

As we are interested in the long time behaviour of solutions, the specific value of the initial data is irrelevant. We therefore 
consider solutions of the Navier--Stokes--Fourier system defined on the open time interval $(T, \infty)$, $T \in R$.

\begin{mdframed}[style=MyFrame]
	
	\begin{Definition} \label{Dw1}
	We say that $(\vr, \vt, \vu)$ is a \emph{weak solution} of the Navier--Stokes--Fourier system \eqref{i1}--\eqref{i3}, 
	\eqref{i7}, \eqref{i8}, with the boundary conditions \eqref{i5}, \eqref{i6} defined on the time interval $(T, \infty)$
	if the following holds:
	
	\begin{itemize}
		
		\item The solution belongs to the following {\bf regularity class}: 
		\begin{align}
			\vr &\in L^\infty(0,T; L^\gamma(\Omega)) \ \mbox{for some}\ \gamma > 1, \br
			\vu &\in L^2_{\rm loc}(T, \infty; W^{1,2}_0 (\Omega; R^3)), \br 
			\vt^{\beta/2} ,\ \log(\vt) &\in L^2(0,T; W^{1,2}(\Omega) \ \mbox{for some}\ \beta \geq 2,\br
			(\vt - \vtB) &\in L^2(0,T; W^{1,2}_0 (\Omega)).
			\label{w6}
			\end{align}
		 
		\item The {\bf equation of continuity} \eqref{i1} along with its renormalization are satisfied in the sense of distributions, 
		specifically,
		\begin{align} 
			\int_T^\infty \intO{ \left[ \vr \partial_t \varphi + \vr \vu \cdot \Grad \varphi \right] } \dt &= 0, 
			\label{w3} \\
			\int_T^\infty \intO{ \left[ b(\vr) \partial_t \varphi + b(\vr) \vu \cdot \Grad \varphi + \Big( 
				b(\vr) - b'(\vr) \vr \Big) \Div \vu \varphi \right] } \dt &=0
			\label{w4}
			\end{align}
for any $\varphi \in C^1_c((T, \infty) \times \Ov{\Omega} )$, and any $b \in C^1(R)$, 	$b' \in C_c(R)$.
\item The {\bf momentum equation} \eqref{i2} is satisfied in the sense of distributions, 
\begin{align}
\int_T^\infty &\intO{ \left[ \vr \vu \cdot \partial_t \bfphi + \vr \vu \otimes \vu : \Grad \bfphi + 
	p \Div \bfphi \right] } \dt \br &= \int_T^\infty \intO{ \left[ \mathbb{S} : \Grad \bfphi - \vr \Grad G \cdot \bfphi \right] } \dt, 
\label{w5}
\end{align}	
for any $\bfphi \in C^1_c((T, \infty) \times \Omega; R^3)$.

\item The internal energy equation \eqref{i3} is replaced by the {\bf entropy inequality}
\begin{align}
	- \int_T^\infty \intO{ \left[ \vr s \partial_t \varphi + \vr s \vu \cdot \Grad \varphi + \frac{\vc{q}}{\vt} \cdot 
		\Grad \varphi \right] } \dt \geq \int_T^\infty \intO{ \frac{\varphi}{\vt} \left[ \mathbb{S} : \Ds \vu - 
		\frac{\vc{q} \cdot \Grad \vt }{\vt} \right] } \dt
	\label{w7} 
	\end{align}
for any $\varphi \in C^1_c((T, \infty) \times \Omega)$, $\varphi \geq 0$; and the {\bf ballistic energy balance}, 
\begin{align}  
	- \int_T^\infty \partial_t \psi	&\intO{ \left[ \frac{1}{2} \vr |\vu|^2 + \vr e - \tvt \vr s \right] } \dt + \int_T^\infty \psi
	\intO{ \frac{\tvt}{\vt}	 \left[ \mathbb{S}: \Ds \vu - \frac{\vc{q} \cdot \Grad \vt }{\vt} \right] } \dt  \br
	&\leq 
	\int_T^\infty \psi \intO{ \left[ \vr \vu \cdot \Grad G - \vr s \vu \cdot \Grad \tvt - \frac{\vc{q}}{\vt} \cdot \Grad \tvt \right] } 
	\label{w8}
\end{align}
for any $\psi \in C^1_c(T; \infty)$, $\psi \geq 0$, and any $\tvt \in C^1([T; \infty) \times \Ov{\Omega})$, 
\[
\tvt > 0,\ \tvt|_{x_3 = 0} = \ThB,\ \tvt|_{x_3 = 1} = \ThU.
\]
		\end{itemize}
	
	\end{Definition}
	
	\end{mdframed}

The existence of global--in--time weak solutions under the constitutive restrictions specified in the forthcoming section was proved in \cite[Theorem 4.2]{ChauFei}. In addition, the weak solutions comply with the weak--strong uniqueness principle and coincide with  strong solutions as soon as they are smooth.

\subsection{Constitutive relations}

Following \cite[Chapters 1,2]{FeNo6A} we consider the equation of state 
\[
p(\vr, \vt) = p_{\rm m} (\vr, \vt) + p_{\rm rad}(\vt), 
\]
where $p_{\rm m}$ is the pressure of a general \emph{monoatomic} gas, 
\begin{equation} \label{con1}
p_{\rm m} (\vr, \vt) = \frac{2}{3} \vr e_{\rm m}(\vr, \vt),
\end{equation}
enhanced by the radiation pressure 
\[
p_{\rm rad}(\vt) = \frac{a}{3} \vt^4,\ a > 0.
\]
Accordingly, the internal energy reads 
\[
e(\vr, \vt) = e_{\rm m}(\vr, \vt) + e_{\rm rad}(\vr, \vt),\ e_{\rm rad}(\vr, \vt) = \frac{a}{\vr} \vt^4.
\]

To identify the specific form of $p_m$ we successively employ several physical principles, see \cite[Chapter 1]{FeNo6A} for details. 

\begin{itemize}

\item {\bf Gibbs' relation} together with \eqref{con1} yield 
\[
p_m (\vr, \vt) = \vt^{\frac{5}{2}} P \left( \frac{\vr}{\vt^{\frac{3}{2}}  } \right)
\]
for a certain $P \in C^1[0,\infty)$.
Consequently, 
\begin{equation} \label{w9}
	p(\vr, \vt) = \vt^{\frac{5}{2}} P \left( \frac{\vr}{\vt^{\frac{3}{2}}  } \right) + \frac{a}{3} \vt^4,\ 
	e(\vr, \vt) = \frac{3}{2} \frac{\vt^{\frac{5}{2}} }{\vr} P \left( \frac{\vr}{\vt^{\frac{3}{2}}  } \right) + \frac{a}{\vr} \vt^4, \ a > 0.
\end{equation}

\item {\bf Hypothesis of thermodynamics stability}, cf. Bechtel, Rooney, and Forrest \cite{BEROFO}, expressed in terms of 
$P$, reads
\begin{equation} \label{w10}
	P(0) = 0,\ P'(Z) > 0 \ \mbox{for}\ Z \geq 0,\ 0 < \frac{ \frac{5}{3} P(Z) - P'(Z) Z }{Z} \leq c \ \mbox{for}\ Z > 0.
\end{equation} 	
In particular, the function $Z \mapsto P(Z)/ Z^{\frac{5}{3}}$ is decreasing, and we suppose 
\begin{equation} \label{w11}
	\lim_{Z \to \infty} \frac{ P(Z) }{Z^{\frac{5}{3}}} = p_\infty > 0.
\end{equation}

\item 
In accordance with Gibbs' relation \eqref{w1}, the associated entropy
takes the form 
\begin{equation} \label{w12}
	s(\vr, \vt) = \mathcal{S} \left( \frac{\vr}{\vt^{\frac{3}{2}} } \right) + \frac{4a}{3} \frac{\vt^3}{\vr}, 
\end{equation}
where 
\begin{equation} \label{w13}
	\mathcal{S}'(Z) = -\frac{3}{2} \frac{ \frac{5}{3} P(Z) - P'(Z) Z }{Z^2} < 0.
\end{equation}
In addition, the {\bf Third law of thermodynamics}, cf. Belgiorno \cite{BEL1}, \cite{BEL2}, requires the entropy to vanish 
as soon as the abolute temperature approaches zero, 
\begin{equation} \label{w14}
\lim_{Z \to \infty} \mathcal{S}(Z) = 0.
\end{equation}

\end{itemize}

Note that \eqref{w10} -- \eqref{w14} imply
\begin{equation} \label{w15}
0 \leq	\vr \mathcal{S} \left( \frac{\vr}{\vt^{\frac{3}{2}}} \right) \leq c \left(1 + \vr \log^+(\vr) + \vr \log^+(\vt) \right).
	\end{equation}

As for the transport coefficients, we suppose that they are continuously differentiable functions satisfying
\begin{align} 
	0 < \underline{\mu}(1 + \vt) &\leq \mu(\vt),\ |\mu'(\vt)| \leq \Ov{\mu}, \br 
	0 \leq \eta (\vt) &\leq \Ov{\eta}(1 + \vt), \br
	0 < \underline{\kappa} (1 + \vt^\beta) &\leq \kappa (\vt) \leq \Ov{\kappa}(1 + \vt^\beta). \label{w16}
	\end{align}
The existence theory developed in \cite{ChauFei} requires 
\begin{equation} \label{w18}
	\beta > 6.
	\end{equation} 

The state equation specified above, together with the fact that the transport coefficients depend on the temperature, 
are pertinent to models of gaseous stars discussed by Bormann \cite{Borm2}, \cite{Borm1}.

\section{Main results}
\label{M}

Our main result states that the Navier--Stokes--Fourier system in the Rayleigh--B\' enard regime \eqref{i5}, \eqref{i6} admits a bounded absorbing set.

\begin{mdframed}[style=MyFrame] 
	\begin{Theorem} [{\bf Bounded absorbing set}] \label{MT1}
		
		Let $\ThB$, $\ThU$ be two strictly positive constants. Let the pressure $p$, the internal energy $e$, and the entropy $s$ 
		satisfy the hypotheses \eqref{w9}--\eqref{w14}. Let the transport coefficients $\mu$, $\eta$, and $\kappa$ satisfy \eqref{w16}, \eqref{w18}. 
		
		Then for any global--in--time weak solution $(\vr, \vt, \vu)$ defined on a time interval $(T, \infty)$, there exists a constant $\mathcal{E}_\infty$ that depends only on $\ThB$, $\ThU$, the total mass of the fluid 
		\[
		M = \intO{ \vr }, 
		\]
		 and the structural properties of 
		$p = p(\vr, \vt)$, $e = a(\vr, \vt)$, $s = s(\vr, \vt)$ such that 
		\begin{equation} \label{MT1}
			{\rm ess} \limsup_{t \to \infty} \intO{ E(\vr, \vt, \vu)(t, \cdot) } \leq \mathcal{E}_\infty. 
			\end{equation}
		
		If, moreover, 
		\[
		{\rm ess} \limsup_{t \to T+} \intO{ E(\vr, \vt, \vu)(t, \cdot) } \leq \mathcal{E}_0 < \infty,
		\]
then the convergence is uniform in $\mathcal{E}_0$. Specifically, for any $\ep > 0$, there exists a time $T(\ep, \mathcal{E}_0)$ such that 
\[
{\rm ess} \sup_{t > T(\ep, \mathcal{E}_0)} \intO{ E(\vr, \vt, \vu) (t, \cdot) } \leq \mathcal{E}_\infty + \ep.		
\]
		\end{Theorem}
	
	\end{mdframed}

\begin{Remark} \label{MR11}
	
	The same result can be shown for a general bounded domain with an arbitrary (nonconstant) profile of the boundary temperature and a general potential volume force $\vc{g} = \Grad G$, $G = G(x)$. In particular, the problem posed in the inclined layer studied e.g. by 
	Daniels et al. \cite{DBPB}
	can be included.
	
	\end{Remark}

The existence of a bounded absorbing set for the isentropic ($p = a \vr^{\gamma}$) Navier--Stokes system with the no--slip
boundary conditions was established in \cite{FP14} under the condition $\gamma > \frac{5}{3}$, see also Wang and Wang \cite{WanWan}. For similar results related to the conservative boundary conditions see \cite{FP9} and the monograph \cite{FeiPr}. The existence of a bounded absorbing set for the Navier--Stokes--Fourier system with general Dirichlet boundary conditions was shown in \cite{FeiKwo} under  rather restrictive assumption postulating a hard sphere equation of state for the pressure. Note that this considerably simplifies the analysis
as \emph{uniform} bounds on the fluid density are {\it a priori} available. It is exactly this missing piece of information that makes the analysis of the present paper much more delicate.

The second result concerns the asymptotic compactness of bounded trajectories. 

\begin{mdframed}[style=MyFrame] 
	
	\begin{Theorem} [{\bf Asymptotic compactness}] \label{MT2} 
		
		\medskip
		
		Under the hypotheses of Theorem \ref{MT1}, let \\ $(\vr_n, \vt_n, \vu_n)_{n=1}^\infty$ be a sequence of weak solutions 
		to the Navier--Stokes--Fourier system in the sense of Definition \ref{Dw1} on the time intervals 
		\[
		(T_n, \infty) ,\ T_n \geq - \infty,\ T_n \to - \infty \ \mbox{as}\ n \to \infty,
		\]
		such that 
		\[
		{\rm ess} \sup_{t \to T_n} \intO{ E(\vr_n, \vt_n, \vu_n )(t, \cdot) } \leq \mathcal{E}_0 ,\ 
		\intO{ \vr } = M > 0,
		\]
		uniformly for $n \to \infty$. 
		
		Then there is a subsequence (not relabelled) such that 
		\begin{align} 
			\vr_n &\to \vr \ \mbox{in}\ C_{\rm weak}([-M, M]; L^{\frac{5}{3}}(\Omega)) \cap C([-M, M] ; L^1(\Omega)), \br
			\vt_n &\to \vt \ \mbox{in}\ L^q((-M,M) ; L^4(\Omega)) \ \mbox{for any}\ 1 \leq q < \infty, \br
			\vu_n &\to \vu \ \mbox{weakly in}\ L^2((-M,M); W^{1,2}(\Omega; R^3)) 
			\label{M1}
			\end{align}
		for any $M > 0$, 
where the limit $(\vr, \vt, \vu)$ is an entire weak solution of the Navier--Stokes--Fourier system defined for $t \in R$ and satisfying 
\begin{equation} \label{M2}
	\intO{ E(\vr, \vt, \vu) (t, \cdot) } \leq \mathcal{E}_{\infty} \ \mbox{for a.a.} \ t \in R.
	\end{equation}		
		
		\end{Theorem}
	
	\end{mdframed}

The heart of the paper is the proof of Theorem \ref{MT1}. Once the uniform bounds on the energy are established, the proof of 
Theorem \ref{MT2} reduces to showing compactness of a sequence of bounded solutions. To certain extent, this is similar to the existence proof, where the only essential issue is the strong (a.e. pointwise) convergence of the densities in \eqref{M1}. Unlike in the existence proof, compactness of the densities at an appropriate ``initial'' time is not available here. Fortunately, this problem is nowadays well understood and we refer the reader to \cite[Section 3, Theorem 3.1]{FanFeiHof} for a detailed proof.

The next section is devoted to the proof of Theorem \ref{MT1}. In view of the hypotheses \eqref{w9}, \eqref{w11}, 
\[
p(\vr, \vt) \approx \vr^{\frac{5}{3}} + \vt^4.
\]
As already pointed out, the exponent $\gamma = \frac{5}{3}$ is critical in the simplified isentropic case. To handle this problem 
we use the fact that {\bf (i)} the gravitational force acting on the fluid is of potential type, and {\bf (ii)} the entropy 
satisfies the Third law of thermodynamics, notably \eqref{w14}.

\section{Dissipativity}
\label{B}

Our goal is to prove Theorem \ref{MT1}. 
Suppose that we are given a global--in--time solution $(\vr, \vt, \vu)$ defined on a time interval $(T, \infty)$. The proof 
of asymptotic boundedness leans on several estimates that follow from the basic physical conservation laws. 
Here and hereafter, we fix $\tvt$ to be the unique solution of the Dirichlet problem 
\begin{equation} \label{B2}
	\Del \tvt = 0 \ \mbox{in}\ \Omega,\ 
	\tvt|_{x_3 = 0} = \ThB,\ \tvt|_{x_3 = 1} = \ThU.
\end{equation}
As $\ThB$, $\ThU$ are constant, we easily compute 
\[
\tvt = \tvt(x_3) = \ThB + x_3 \left( {\ThU} - \ThB \right).
\]
Obviously, the same ansatz can be used in the case of general $x-$dependent boundary data.

\subsection{Mass conservation}

It follows from the equation of continuity \eqref{w3} that the total mass of the fluid is a constant of motion, 
\begin{equation} \label{B2b}
	M = \intO{ \vr (t, \cdot) } \ \mbox{for any}\ t > T.
\end{equation}

In addition, as the volume force is potential, we can write 
\[
\intO{ \vr \vu \cdot \Grad G } = \frac{\D }{\dt} \intO{ \vr G },\ G = - x_3.
\]
Consequently, the ballistic energy balance \eqref{w8} takes the form 
\begin{align}  
	\frac{\D }{\dt}	&\intO{ \left[ \frac{1}{2} \vr |\vu|^2 + \vr e - \tvt \vr s - \vr G \right] }  + 
	\intO{ \frac{\tvt}{\vt}	 \left[ \mathbb{S}: \Ds \vu - \frac{\vc{q} \cdot \Grad \vt }{\vt} \right] }   \br
	&\leq 
	- \intO{ \left[ \vr s \vu \cdot \Grad \tvt + \frac{\vc{q}}{\vt} \cdot \Grad \tvt \right] } \ \mbox{in}\ 
	\mathcal{D}'(T, \infty). 
	\label{B1}
\end{align}
It is worth--noting that the same argument applies for a general Lipschitz potential $G = G(x)$.

\subsection{Coercivity of the dissipative term}

It follows from the hypotheses \eqref{w16} and Korn--Poincar\' e inequality that
\begin{align}
&\intO{ \frac{\tvt}{\vt}	 \left( \mathbb{S}(\vt, \Ds \vu) : \Ds \vu - \frac{\vc{q}(\vt, \Grad \vt) \cdot \Grad \vt }{\vt} \right) }
\br &\quad \geq c \inf \{ \ThU, \ThB \} \left( \| \vu \|^2_{W^{1,2}(\Omega; R^3)} + \| \Grad \vt^{\frac{\beta}{2}} \|^2_{L^2(\Omega; R^3)} 
+ \| \Grad \log(\vt) \|^2_{L^2(\Omega; R^3)}  \right).
\nonumber
\end{align}
Consequently, adding the boundary integrals to the left--hand side and using Poincar\' e inequality, we get 
\begin{align}
&\left( \| \vu \|^2_{W^{1,2}(\Omega; R^3)} + \| \vt^{\frac{\beta}{2}} \|^2_{W^{1,2}(\Omega)} 
+ \| \log(\vt) \|^2_{W^{1,2}(\Omega)}  \right) \br	
	&\quad \leq  c( \ThU, \ThB  ) \left[ 1 + \intO{ \frac{\tvt}{\vt}	 \left( \mathbb{S}(\vt, \Ds \vu) : \Ds \vu - \frac{\vc{q}(\vt, \Grad \vt) \cdot \Grad \vt }{\vt} \right) } \right]
		\label{B2a}
\end{align}

\subsection{Energy estimates}

To simplify the ballistic energy inequality \eqref{B1}, we first realize, by virtue of \eqref{B2}, 
\[
\intO{ \frac{\vc{q} (\vt, \Grad \vt)}{\vt} \cdot \Grad \tvt } = - 
\intO{ \frac{\kappa(\vt)}{\vt} \Grad \vt \cdot \Grad \tvt } = \int_{\partial \Omega} \mathcal{K} (\tvt) \Grad \tvt \ \D \sigma_x,
\]
where 
\[
\mathcal{K}' (\vt) = \kappa(\vt).
\] 
Consequently, the ballistic energy inequality \eqref{B1} reduces to 
	\begin{align}  
	\frac{\D }{\dt}	&\intO{ \left( \frac{1}{2} \vr |\vu|^2 + \vr e - \tvt \vr s - \vr G \right) } + 
	\intO{ \frac{\tvt}{\vt}	 \left( \mathbb{S}(\vt, \Ds \vu) : \Ds \vu - \frac{\vc{q}(\vt, \Grad \vt) \cdot \Grad \vt }{\vt} \right) }  \br
	&\leq 
	- \intO{ \left[ \vr s \vu \cdot \Grad \tvt  \right] } + c( \ThU, \ThB  ) .
	\label{B1bisa}
\end{align}

\subsection{Entropy estimates}

In accordance with hypothesis \eqref{w12}, 
\[
\intO{ \vr s \vu \cdot \nabla\tvt } = \intO{ \vr \mathcal{S} \left( \frac{\vr}{\vt^{\frac{3}{2}}} \right) \vu \cdot \Grad \tvt } + 
\frac{4a}{3} \intO{ \vt^3 \vu \cdot \nabla\tvt }, 
\]
where, by virtue of \eqref{w13}, \eqref{w14},
\begin{equation} \label{BB1}
\vr \mathcal{S} \left( \frac{\vr}{\vt^{\frac{3}{2}}} \right) \leq \vr \mathcal{S}(r) \ \mbox{provided} \ \frac{\vr}{\vt^{\frac{3}{2}}} 
\geq r,\ \mbox{where}\ \mathcal{S}(r) \to 0 \ \mbox{as}\ r \to \infty.
\end{equation}
If 
\[
\frac{\vr}{\vt^{\frac{3}{2}}} < r \ \mbox{or}\ \vr < r \vt^{\frac{3}{2}}, 
\]
we get, by virtue of \eqref{w15}, 
\begin{equation} \label{BB2}
0 \leq \vr \mathcal{S} \left( \frac{\vr}{\vt^{\frac{3}{2}}} \right) 
\leq c \left(1 + r \vt^{\frac{3}{2}} \left[ \log^+ (r \vt^{\frac{3}{2}} ) + \log^+(\vt) \right] \right). 
\end{equation}
As $\beta > 6$ in hypothesis \eqref{w16}, we may combine \eqref{B2a} with \eqref{B1bisa}--\eqref{BB2} to obtain  
	\begin{align}  
	\frac{\D }{\dt}	&\intO{ \left( \frac{1}{2} \vr |\vu|^2 + \vr e - \tvt \vr s  - \vr G \right) } \br &+
	c_1 ( \ThU, \ThB  )  \left( \| \vu \|^2_{W^{1,2}(\Omega; R^3)} + \| \vt^{\frac{\beta}{2}} \|^2_{W^{1,2}(\Omega)} 
	+ \| \log(\vt) \|^2_{W^{1,2}(\Omega)}  \right) 
	 \br
	&\leq c_2 ( \ThU, \ThB  )  \mathcal{S}(r)  
	\intO{ \vr |\vu| } + \Lambda(\ThU, \ThB  , r ) ,\br & \quad \quad \mbox{where}\ c_1 > 0,\ c_2 > 0,\ \mbox{and}\ \Lambda(r) \to \infty 
	\ \mbox{if}\ r \to \infty.
	\label{B2bis}
\end{align}

The main problem to conclude is the fact that forcing term $\intO{ \vr |\vu| }$ on the right--hand side is not directly controlled by the dissipation on the left--hand side. To this end, we need the so--called pressure estimates which we recall in the next section.

\subsection{Pressure estimates}
\label{Bb}

To continue, we recall the inverse of the divergence known as Bogovskii operator:
\begin{align}
	\mathcal{B} : L^q_0 (\Omega) &\equiv \left\{ f \in L^q(\Omega) \ \Big| \ \intO{ f } = 0 \right\}
	\to W^{1,q}_0(\Omega, R^d),\ 1 < q < \infty,\br 
	\Div \mathcal{B}[f] &= f, \br 
	\| \mathcal{B}[ \Div \vc{g} ] \|_{L^r(\Omega)} &\leq c \| \vc{g} \|_{L^r(\Omega)},\ 1 < r < \infty 
	\ \mbox{whenever}\ \vc{g} \cdot \vc{n}|_{\partial \Omega} = 0.
	\label{Bog}
\end{align}
see e.g. Galdi \cite[Chapter 3]{GALN} or Gei{\ss}ert, Heck, and Hieber\cite{GEHEHI}. 

Now, the test function 
\[
\varphi (t,x) = \mathcal{B} \left[ b(\vr) - \frac{1}{|\Omega|} \intO{ b(\vr) } \right]
\]
in the momentum equation yields
\begin{align}
	\int_\tau^{\tau + 1} &\intO{ p(\vr, \vt) b(\vr) } \dt = \int_\tau^{\tau + 1} \frac{1}{|\Omega|} \left( \intO{ b(\vr) }  \right) \left( \intO{ p(\vr, \vt)  } \right) \dt \br 
	&-\int_\tau^{\tau + 1} \intO{ \vr (\vu \otimes \vu): \Grad \mathcal{B} \left[ b(\vr) - \frac{1}{|\Omega|} \intO{ b(\vr) } \right] } 
	\dt \br 
	&+ \int_\tau^{\tau + 1} \intO{ \mathbb{S}(\vt, \Ds \vu) : \Grad \mathcal{B} \left[ b(\vr) - \frac{1}{|\Omega|} \intO{ b(\vr) } \right] } 
	\dt\br & - \int_\tau^{\tau + 1} \intO{ \vr \Grad G \cdot \mathcal{B} \left[ b(\vr) - \frac{1}{|\Omega|} \intO{ b(\vr) } \right] } \dt \br 
	&+ \left[ \intO{ \vr \vu \cdot \mathcal{B} \left[ b(\vr) - \frac{1}{|\Omega|} \intO{ b(\vr) } \right] } \right]_{t = \tau}^{t = \tau + 1} \br 
	&- \int_\tau^{\tau + 1} \intO{  \vr \vu \cdot \partial_t \mathcal{B} \left[ b(\vr) - \frac{1}{|\Omega|} \intO{ b(\vr) } \right] } 
	\dt.
	\label{B3}	
\end{align}
In addition, as $\vr$ satisfies the renormalized equation of continuity, we obtain
\begin{align}
\int_\tau^{\tau + 1} &\intO{ \vr \vu \cdot \partial_t \mathcal{B} \left[ b(\vr) - \frac{1}{|\Omega|} \intO{ b(\vr) } \right] } 
\dt \br & = - \int_\tau^{\tau + 1} \intO{	\vr \vu \cdot \mathcal{B}[\Div (b(\vr) \vu ] } \dt \br 
&+ \int_{\tau}^{\tau + 1} \intO{ \vr \vu \cdot \mathcal{B} \left[ (b(\vr) - b'(\vr) \vr) \Div \vu - 
	\frac{1}{|\Omega|} \intO{ (b(\vr) - b'(\vr) \vr) \Div \vu } \right] } \dt.
\label{B4}	
	\end{align}
The unit length of the time interval has been chosen just for simplicity.

\subsection{Uniform bounds}

In view of the structural restrictions imposed through hypotheses \eqref{w10}, \eqref{w12} and \eqref{w15}, for any $\lambda > 1$ there exist 
two constants $c_1(\lambda)$, $c_2(\lambda)$ such that  
\begin{equation} \label{B4a}
	c_1(\lambda, {\rm data}) + \frac{1}{\lambda} E(\vr, \vt, \vu) \leq E_{\tvt}(\vr, \vt, \vu) - \vr G \leq \lambda E(\vr, \vt, \vu) 
	+ c_2(\lambda, {\rm data}).
	\end{equation} 

The following result is crucial for showing the existence of a bounded absorbing set. 

\begin{Lemma} \label{LB1}
	
	Suppose that 
	\begin{equation} \label{B7}
	\intO{ \left[ E_{\tvt} (\tau, \cdot) - \vr(\tau, \cdot) G \right] } - \intO{ \left[ E_{\tvt}(\tau
		 + 1, \cdot) - \vr(\tau + 1, \cdot) G \right] } \leq K.
	\end{equation}

Then there exists $L = L(K, M, {\rm data})$ such that 
\begin{equation} \label{B8}
{\rm ess} \sup_{\tau \leq t \leq \tau + 1} \intO{ E(t, \cdot) } \leq L.
\end{equation}	
	
	\end{Lemma}

\begin{Remark} \label{RB1} 
	Strictly speaking, the pointwise values of the ballistic energy appearing in \eqref{B7} are defined only for a.a. $\tau \in (T; \infty)$. However, thanks to the inequality \eqref{B2bis}, we may identify $E_{\tvt}$ with its, say, c\` agl\` ad representative defined for any $\tau > T$.  
	
	\end{Remark}

The rest of this subsection is devoted to the proof of Lemma \ref{LB1}.
If \eqref{B7} holds, it follows from \eqref{B2bis} that
\begin{align}  
	\int_{\tau}^{\tau + 1}& 
	\left( \| \vu \|^2_{W^{1,2}(\Omega; R^3)} + \| \vt^{\frac{\beta}{2}} \|_{W^{1,2}(\Omega)} + \| \log(\vt) \|^2_{W^{1,2}(\Omega)} \right) \dt \br
	&\leq \Cdk \left( 1 + \mathcal{S}(r)
	\int_{\tau}^{\tau + 1} \intO{ \vr  |\vu|   } \dt \right) + \Lambda ({\rm data}, K, r ).
	\label{B10}
\end{align}

\subsubsection{Pressure estimates revisited}

At this stage, we use the pressure estimates \eqref{B3}, with 
\[
b(\vr) = \vr^{\alpha},\ \alpha > 0.
\]
In view of the hypotheses \eqref{w9}, \eqref{w11} imposed on the equation of state, we have 
\begin{equation} \label{B13}
c_1 \Big( \vr^{\frac{5}{3}} + \vt^4 \Big) \leq p(\vr, \vt) \leq c_2 \Big(\vr^{\frac{5}{3}} + \vt^4 + 1 \Big), \ c_1, c_2 > 0. 	
	\end{equation}
Moreover, as the total mass is constant via \eqref{B2b}, the smoothing properties of $\mathcal{B}$ stated in \eqref{Bog} imply 
\[
\left| \mathcal{B} \left[ \vr^\alpha - \frac{1}{|\Omega|} \intO{ \vr^\alpha  } \right] \right| \leq c(M)
\ \mbox{as soon as}\ \alpha < \frac{1}{3} 
\]
to provide that $W^{1,\frac{1}{\alpha} }(\Omega)\subset L^\infty(\Omega)$. Thus inequality \eqref{B3} gives rise to 
	
\begin{align}
	\int_\tau^{\tau + 1} &\intO{ \vr^{\frac{5}{3} + \alpha} } \dt \leq c(M) \left( 1 + \int_\tau^{\tau + 1} \intO{ \vt^4 }  \dt \right. \br 
	&-\int_\tau^{\tau + 1} \intO{ \vr (\vu \otimes \vu): \Grad \mathcal{B} \left[ \vr^\alpha - \frac{1}{|\Omega|} \intO{ \vr^\alpha } \right] } 
	\dt \br 
	&+ \int_\tau^{\tau + 1} \intO{ \mathbb{S}(\vt, \Ds \vu) : \Grad \mathcal{B} \left[ \vr^\alpha - \frac{1}{|\Omega|} \intO{ \vr^\alpha } \right] } 
	\dt\br 
	&+ \left[ \intO{ \vr \vu \cdot \mathcal{B} \left[ \vr^\alpha - \frac{1}{|\Omega|} \intO{ \vr^\alpha } \right] } \right]_{t = \tau}^{t = \tau + 1} \br 
	&- \left. \int_\tau^{\tau + 1} \intO{ \vr \vu \cdot \partial_t \mathcal{B} \left[ \vr^\alpha - \frac{1}{|\Omega|} \intO{ \vr^\alpha } \right] } 
	\dt \right) .
	\label{B14}	
\end{align}

Next, using again the smoothing properties \eqref{Bog} of $\mathcal{B}$ we get
\begin{align}
&\left| \int_\tau^{\tau + 1} \intO{ \vr (\vu \otimes \vu): \Grad \mathcal{B} \left[ \vr^\alpha - \frac{1}{|\Omega|} \intO{ \vr^\alpha } \right] } \dt \right| \br
&\quad \leq \int_{\tau}^{\tau + 1} \| \vr \|_{L^\gamma(\Omega)} \| \vu \|^2_{L^6(\Omega; R^3)} \| \vr^\alpha \|_{L^q(\Omega)}  
\dt	\br &\quad \leq \sup_{t \in (\tau, \tau + 1)} \| \vr \|_{L^\gamma (\Omega)} \int_{\tau}^{\tau + T} \| \vu \|^2_{W^{1,2}(\Omega; R^3)} \sup_{t \in (\tau, \tau + 1)} \| \vr^\alpha \|_{L^q (\Omega)}
\dt,
	\label{B15}
	\end{align}
where 
\[
q = \frac{3 \gamma}{2 \gamma - 3} > 1 \ \mbox{provided}\ \gamma > \frac{3}{2}. 
\]
Thus setting 
\begin{equation} \label{B16}
\gamma = \frac{5}{3},\ 	
\alpha = \frac{2 \gamma - 3}{3 \gamma} = \frac{1}{15} < \frac{1}{3},
\end{equation}
we may use the total mass conservation \eqref{B2b} to conclude
\begin{align}
	&\left| \int_\tau^{\tau + 1} \intO{ \vr (\vu \otimes \vu): \Grad \mathcal{B} \left[ \vr^\alpha - \frac{1}{|\Omega|} \intO{ \vr^\alpha } \right] } \dt \right| 	\br &\quad \leq c(M) \sup_{t \in (\tau, \tau + 1)} \| \vr \|_{L^{\frac{5}{3}} (\Omega)} \int_{\tau}^{\tau + 1} \| \vu \|^2_{W^{1,2}(\Omega; R^3)} 
	\dt,
	\label{B17}
\end{align}

Similarly, going back to \eqref{B4} we have
\begin{align}
\left| \int_{\tau}^{\tau + 1} \intO{ \vr \vu \cdot \mathcal{B}[\Div (\vr^\alpha \vu)] } \right|	 \leq 
\int_0^{\tau + 1} \| \vr \|_{L^\gamma (\Omega)} \| \vu \|_{L^6(\Omega; R^3)} \| \vr^\alpha \vu \|_{L^q(\Omega; R^3)} \dt,
\nonumber
	\end{align}
where
\[
\frac{1}{\gamma}+ \frac{1}{6} + \frac{1}{q} = 1.
\]
Moreover, 
\[
\| \vr^\alpha \vu \|_{L^q(\Omega; R^3)} \leq \| \vu \|_{L^6(\Omega; R^3)} \| \vr^\alpha \|_{L^p(\Omega)}, 
\mbox{where}\ 
\frac{1}{q} = \frac{1}{6} + \frac{1}{p};  
\]
whence 
\begin{align}
	\left| \int_{\tau}^{\tau + 1} \intO{ \vr \vu \cdot \mathcal{B}[\Div (\vr^\alpha \vu)] } \right|	 \leq 
c(M) \sup_{t \in (\tau, \tau + 1)} \| \vr \|_{L^{\frac{5}{3}} (\Omega)} \int_{\tau}^{\tau + 1} \| \vu \|^2_{W^{1,2}(\Omega; R^3)} \dt	
	\label{B18}
\end{align}
as soon as  \eqref{B16} holds.

Finally, 
\begin{align}
&\left|  \int_{\tau}^{\tau + 1} \intO{ \vr \vu \cdot \mathcal{B} \left[ \vr^\alpha \Div \vu - 
	\frac{1}{|\Omega|} \intO{ \vr^\alpha \Div \vu } \right] } \dt \right| \br &\quad \leq 
\int_\tau^{\tau + 1} \| \vr \|_{L^\gamma (\Omega)} \| \vu \|_{L^6(\Omega; R^3)} \left\| \mathcal{B} \left[ \vr^\alpha \Div \vu - 
\frac{1}{|\Omega|} \intO{ \vr^\alpha \Div \vu } \right] \right\|_{L^q(\Omega; R^3)} \dt,
\nonumber
\end{align}
where
\[
\frac{1}{\gamma} + \frac{1}{6} + \frac{1}{q} = 1.
\]
Furthermore, 
\[
\left\| \mathcal{B} \left[ \vr^\alpha \Div \vu - 
\frac{1}{|\Omega|} \intO{ \vr^\alpha \Div \vu } \right] \right\|_{L^q(\Omega; R^3)} \aleq 
\| \vr^\alpha \Div \vu \|_{L^r(\Omega; R^3)},\ q = \frac{3 r}{3 - r}, 
\]
and 
\[
\| \vr^\alpha \Div \vu \|_{L^r(\Omega; R^3)} \leq \| \vu \|_{W^{1,2}(\Omega; R^3)} \| \vr^\alpha \|_{L^p(\Omega)},\ 
\mbox{with}\ \frac{1}{2} + \frac{1}{p} = \frac{1}{r}.
\]
Consequently, condition \eqref{B16} yields 
\begin{align}
	&\left|  \int_{\tau}^{\tau + 1} \intO{ \vr \vu \cdot \mathcal{B} \left[ \vr^\alpha \Div \vu - 
		\frac{1}{|\Omega|} \intO{ \vr^\alpha \Div \vu } \right] } \dt \right| \br &\quad \leq 
	c(M) \sup_{t \in (\tau, \tau + 1)} \| \vr \|_{L^{\frac{5}{3}} (\Omega)} \int_{\tau}^{\tau + 1} \| \vu \|^2_{W^{1,2}(\Omega; R^3)} \dt.
	\label{B19}
\end{align}

Summing up the previous inequalities and going back to \eqref{B14}, we get 
\begin{align}
	\int_\tau^{\tau + 1} &\intO{ \vr^{\frac{5}{3} + \alpha} } \dt \leq c(M) \left( 1 + \int_\tau^{\tau + 1} \intO{ \vt^4 }  \dt \right. \br 
	&+	\sup_{t \in (\tau, \tau + 1)} \| \vr \|_{L^{\frac{5}{3}} (\Omega)} \int_{\tau}^{\tau + 1} \| \vu \|^2_{W^{1,2}(\Omega; R^3)}
	\dt
	\br 
	&+ \int_\tau^{\tau + T} \intO{ \mathbb{S}(\vt, \Ds \vu) : \Grad \mathcal{B} \left[ \vr^\alpha - \frac{1}{|\Omega|} \intO{ \vr^\alpha } \right] } 
	\dt\br 
	&+ \left. \left[ \intO{ \vr \vu \cdot \mathcal{B} \left[ \vr^\alpha - \frac{1}{|\Omega|} \intO{ \vr^\alpha } \right] } \right]_{t = \tau}^{t = \tau + 1} \right),\ \alpha = \frac{1}{15}.
	\label{B20}	
\end{align}	

Now, 
\begin{align} 
&\intO{ \mathbb{S} (\vt, \Ds \vu) : \Grad \mathcal{B} \left[ \vr^\alpha - \frac{1}{|\Omega|} \intO{ \vr^\alpha } \right] } \br 
&\quad \leq  (1 + \| \vt \|_{L^4(\Omega)} ) \| \vu \|_{W^{1,2}(\Omega; R^3)} \left\| 
 \Grad \mathcal{B} \left[ \vr^\alpha - \frac{1}{|\Omega|} \intO{ \vr^\alpha } \right] \right\|_{L^4(\Omega; R^3)} \br 
&\quad \leq c(M) (1 + \| \vt \|_{L^4(\Omega)} ) \| \vu \|_{W^{1,2}(\Omega; R^3)}.
\nonumber 
\end{align} 
We therefore conclude 
\begin{align}
	\int_\tau^{\tau + 1} &\intO{ \vr^{\frac{5}{3} + \alpha} } \dt \leq c(M) \left[ 1 + \int_\tau^{\tau + 1} \intO{ \vt^4 }  \dt \right. \br 
	&+	\left(1 + \sup_{t \in (\tau, \tau + 1)} \| \vr \|_{L^{\frac{5}{3}} (\Omega)} \right) \int_{\tau}^{\tau + 1} \| \vu \|^2_{W^{1,2}(\Omega; R^3)}
	\dt
	+ \left. \sup_{t \in (\tau, \tau + 1)} \intO{ \vr |\vu| } \right],\ \alpha = \frac{1}{15}. 
	\label{B21}	
\end{align}

\subsubsection{Proof of Lemma \ref{LB1}}

Now, in accordance with \eqref{B10},
\begin{align}
 \int_\tau^{\tau + 1} \intO{ \vt^4 } &\leq \Cd \left(1 +  \int_\tau^{\tau + 1} \| \vt^{\frac{\beta}{2}} \|^2_{W^{1,2}(\Omega)} 
 \dt \right) \br &\leq c({\rm data},K) \left(1 + \mathcal{S}(r) \int_{\tau}^{\tau + 1} \intO{ \vr |\vu| } \dt \right) + 
 \Lambda({\rm data}, K, r),
 \nonumber	
	\end{align}
where we may fix $r = 1$.
Consequently, inequality \eqref{B21} reduces to 
\begin{align}
	\int_\tau^{\tau + 1} \intO{ \vr^{\frac{5}{3} + \alpha} } \dt &\leq c(K,M, {\rm data}) \left[ \left(1 + \sup_{t \in (\tau, \tau + 1)} \| \vr \|_{L^\frac{5}{3} (\Omega)} \right) \int_{\tau}^{\tau + 1} \| \vu \|^2_{W^{1,2}(\Omega; R^3)}
	\dt \right. \br 
	&+  \left. \sup_{t \in (\tau, \tau + 1)} \intO{ \vr |\vu| } + 1 \right],\ \alpha = \frac{1}{15}.
	\label{B23}	
\end{align}

Next, it follows from the hypotheses \eqref{B7}, \eqref{B10}, and \eqref{B2bis} that 
\begin{equation} \label{B25}
	\sup_{t \in (\tau, \tau + 1)} \intO{ E(t, \cdot) } \leq \Cd \left( 1 + \int_{\tau}^{\tau + 1} E(s, \cdot) \D s \right).
\end{equation}
Moreover, relation \eqref{B10} yields 
\[
\int_{\tau}^{\tau+1} \| \vu \|^2_{W^{1,2}(\Omega; R^3)} \leq c({\rm data}, K) \mathcal{S}(r) \int_{\tau}^{\tau + 1} \intO{ \vr |\vu| } \dt 
+ \Lambda ({\rm data},K, r),
\]
where, by H\" older's inequality and Sobolev embedding theorem, 
\[
\intO{ \vr |\vu| } \leq \| \sqrt{\vr} \|_{L^2(\Omega)} \| \sqrt{\vr} \|_{L^3(\Omega)} \| \vu \|_{L^6(\Omega; R^3)} 
\leq c \sqrt{M} \| \sqrt{\vr} \|_{L^3(\Omega)} \| \vu \|_{W^{1,2}(\Omega; R^3)}.
\]
Thus we may infer that 
\begin{equation} \label{B26}
\int_{\tau}^{\tau+1} \| \vu \|^2_{W^{1,2}(\Omega; R^3)} \leq c({\rm data}, K,M) \mathcal{S}(r) \int_{\tau}^{\tau + 1} \| \vr \|_{L^{\frac{3}{2}}(\Omega)}
+ \Lambda ({\rm data},K, r).
\end{equation}
Finally, using \eqref{B26} we may estimate the kinetic energy,
\begin{align}
\int_{\tau}^{\tau + 1} & \intO{ \vr |\vu|^2 } \dt \leq \sup_{t \in (\tau, \tau + 1)} \| \vr \|_{L^{\frac{3}{2}}(\Omega)} 
\int_{\tau}^{\tau + 1} \| \vu \|^2_{L^6(\Omega; R^3)} \dt \br &\leq c \sup_{t \in (\tau, \tau + 1)} \| \vr \|_{L^{\frac{3}{2}}(\Omega)} 
\int_{\tau}^{\tau + 1} \| \vu \|^2_{W^{1,2}(\Omega; R^3)} \dt \br 
&\leq \Lambda ({\rm data},K, r) \sup_{t \in (\tau, \tau + 1)} \| \vr \|_{L^{\frac{3}{2}}(\Omega)} + 
c({\rm data}, K,M) \mathcal{S}(r) \sup_{t \in (\tau, \tau + 1)} \| \vr \|_{L^{\frac{3}{2}}(\Omega)}\int_{\tau}^{\tau + 1} \| \vr \|_{L^{\frac{3}{2}}(\Omega)} \dt. \nonumber
\end{align}
Now, by interpolation, 
\[
\| \vr \|_{L^{\frac{3}{2}}(\Omega)} \leq \| \vr \|_{L^{\frac{5}{3}}(\Omega)}^{\frac{5}{6}} \| \vr \|_{L^1(\Omega)}^{\frac{1}{6}};
\]	
whence 
\begin{align} 
&\int_{\tau}^{\tau + 1}  \intO{ \vr |\vu|^2 }  \br 
&\quad \leq \Lambda ({\rm data},K, M, r) \sup_{t \in (\tau, \tau + 1)} \| \vr \|_{L^{\frac{5}{3}}(\Omega)}^{\frac{5}{6}} 
+ c({\rm data}, K,M) \mathcal{S}(r) \sup_{t \in (\tau, \tau + 1)} \| \vr \|_{L^{\frac{5}{3}}(\Omega)}^{\frac{5}{6}} \int_{\tau}^{\tau + 1} \| \vr \|_{L^{\frac{5}{3}}(\Omega)}^{\frac{5}{6}} \dt.
\label{B28}
\end{align}	

Going back to \eqref{B23} and using \eqref{B26} we get
\begin{align}
	\int_\tau^{\tau + 1} \intO{ \vr^{\frac{5}{3} + \alpha} } \dt &\leq \Lambda(K,M, {\rm data},r) \left[ \left(1 + \sup_{t \in (\tau, \tau + 1)} \| \vr \|_{L^\frac{5}{3} (\Omega)} \right)\dt \right. \br 
	&+ c(K,M, {\rm data}) \mathcal{S}(r) \left(1 + \sup_{t \in (\tau, \tau + 1)} \| \vr \|_{L^\frac{5}{3} (\Omega)} \right)
		\int_\tau^{\tau + 1} \| \vr \|_{L^{\frac{3}{2}}(\Omega)} \dt
	 \br 
	&+  \left. \sup_{t \in (\tau, \tau + 1)} \intO{ \vr |\vu| } + 1 \right] \br 
	 &\leq \Lambda(K,M, {\rm data},r) \left[ \left(1 + \sup_{t \in (\tau, \tau + 1)} \| \vr \|_{L^\frac{5}{3} (\Omega)} \right)\dt \right. \br 
	&+ c(K,M, {\rm data}) \mathcal{S}(r) \left(1 + \sup_{t \in (\tau, \tau + 1)} \| \vr \|_{L^\frac{5}{3} (\Omega)} \right)
	\int_\tau^{\tau + 1} \| \vr \|_{L^{\frac{5}{3}}(\Omega)}^{\frac{5}{6}} \dt
	\br 
	&+  \left. \sup_{t \in (\tau, \tau + 1)} \intO{ \vr |\vu| } + 1 \right],\ \alpha = \frac{1}{15}.
	\label{B29}	
\end{align}
Now, interpolating $L^1$ and $L^{\frac{5}{3} + \alpha}$, we get 
\[
\int_{\tau}^{\tau + 1} \intO{ \vr^{\frac{5}{3}} } \dt \leq c(M) \left( \int_{\tau}^{\tau + 1} \intO{ \vr^{\frac{5}{3} + \alpha } }
\dt \right)^{\frac{10}{11}} \ \mbox{provided}\ \alpha = \frac{1}{15}.
\]

Gathering the available bounds we conclude 
\begin{align}
\sup_{t \in (\tau, \tau + 1)} &\intO{ E(t, \cdot) } \leq \Cd \left( 1 + \int_{\tau}^{\tau + 1} E(s, \cdot) \D s \right) \br 
\leq  \Cd &\left( 1 + \int_{\tau}^{\tau + 1} 
\left( \| \vu \|^2_{W^{1,2}(\Omega; R^3)} + \| \vt^{\frac{\beta}{2}} \|_{W^{1,2}(\Omega)} + \| \log(\vt) \|^2_{W^{1,2}(\Omega)} \right) \dt \right)	\br 
+ \Cd& \left( \int_{\tau}^{\tau + 1} \intO{ \vr |\vu|^2 } \dt + \int_{\tau}^{\tau + 1} \intO{ \vr^{\frac{5}{3}} } \dt \right) \br 
\leq \Lambda({\rm data}, K,r)& \left[ 1 + \left( \sup_{t \in (\tau, \tau + 1)} \intO{ E } \dt \right)^{\lambda} \right] 
+ c({\rm data}, K, M) \mathcal{S}(r)  \sup_{t \in (\tau, \tau + 1)} \intO{ E } \dt
\label{B30}
	\end{align}
for certain $0 < \lambda < 1$.
Consequently, choosing $r = r({\rm data}, K, M)$ large enough, the desired conclusion follows since $\mathcal{S}(r) \to 0$ 
as $r \to \infty$.

We have proved Lemma \ref{LB1}.

\subsection{Bounded absorbing sets}

The existence of a bounded absorbing set follows easily from Lemma \ref{LB1}. Indeed consider a global--in--time solution as in 
Theorem \ref{MT1} satisfying 
\[
{\rm ess} \limsup_{t \to T-} \intO{ E(\vr, \vt, \vu)(t,\cdot) } \leq \mathcal{E}_0.
\]
Consider $K = 1$ in Lemma \ref{LB1}. In view of \eqref{B4a}, there exists $\tau = \tau(\mathcal{E}_0)$ such that \eqref{B4a} 
holds. Repeating the same argument with $\mathcal{E}_0$ replaced by $L$ given by \eqref{B8} we deduce that there exists $H = H(L)$ such that any time interval $(s, s + H)$, $s \geq \tau(\mathcal{E}_0)$ contains $\tau$ such that \eqref{B7} holds yielding
\eqref{B8}. Consequently, we may consider 
\[
\mathcal{E}_\infty \leq L \exp H(L). 
\] 

We have proved Theorem \ref{MT1}. As pointed out in Section \ref{M}, Theorem \ref{MT1} yields Theorem \ref{MT2} via the existing 
arguments presented e.g. in \cite{FanFeiHof}.

\section{Applications, long--time behavior}
\label{L}

We finish the paper by discussing the impact of Theorems \ref{MT1}, \ref{MT2} on the long time behavior of solutions to the 
Rayleigh--B\' enard problem in the framework of compressible viscous and heat conducting fluids.

\subsection{Trajectory space} 

We start by introducing a suitable \emph{trajectory space} $\mathcal{T}$. In view of the framework of Theorems \ref{MT1}, \ref{MT2}, the ``natural'' trajectory space should be based on the standard phase variables $(\vr, \vt, \vu)$. Unfortunately, neither $\vt$ nor $\vu$ admit well defined \emph{instantaneous values} at any time $t \in R$. It is therefore more convenient to consider the \emph{conservative entropy} variables $(\vr, S, \vm)$, with 
\[
\mbox{momentum}\ \vm = \vr \vu,\ \mbox{and total entropy}\ S = \vr s(\vr, \vt).
\]
On the one hand, the state variables $(\vr, S, \vm)$ are uniquely determined by $(\vr, \vt, \vu)$. On the other hand, knowing 
$(\vr, S, \vm)$ we first obtain $\vt$ as $\vt \mapsto s(\vr, \vt)$ is a strictly increasing function. The velocity $\vu$ 
is {\it a priori} not well defined on the hypothetical vacuum zone, however, it can be recovered in terms of $(\vr, \vt, \vm)$ 
from the momentum equation  \eqref{w5}. 

The phase variables $(\vr, S, \vm)$ admit well defined instantaneous values understood in the weak sense. Specifically, 
it follows from the weak formulation \eqref{w3}, \eqref{w5}, \eqref{w7} that the one sided limits
\begin{align}
\left< \vr(\tau-, \cdot); \phi \right> &\equiv \lim_{\delta \to 0+} \frac{1}{\delta} \int_{\tau - \delta}^\tau \intO{ \vr(t,\cdot) \phi } \dt, \ 
\left< \vr(\tau+, \cdot); \phi \right> \equiv \lim_{\delta \to 0+} \frac{1}{\delta} \int_{\tau}^{\tau + \delta} \intO{ \vr(t,\cdot) \phi } \dt \br
\left< \vm(\tau-, \cdot); \bfphi \right> &\equiv \lim_{\delta \to 0+} \frac{1}{\delta} \int_{\tau - \delta}^\tau \intO{ \vm(t,\cdot) \cdot \bfphi } \dt, \ 
\left< \vm(\tau+, \cdot) ; \bfphi \right> \equiv \lim_{\delta \to 0+} \frac{1}{\delta} \int_{\tau}^{\tau + \delta} \intO{ \vm(t,\cdot) \cdot \bfphi } \dt, \br
\left< S(\tau-, \cdot); \phi \right> &\equiv \lim_{\delta \to 0+} \frac{1}{\delta} \int_{\tau - \delta}^\tau \intO{ \vr s (t,\cdot) \phi } \dt, \ 
\left< S(\tau+, \cdot); \phi \right> \equiv \lim_{\delta \to 0+} \frac{1}{\delta} \int_{\tau}^{\tau + \delta} \intO{ \vr s(t,\cdot) \phi } \dt	
\nonumber
	\end{align}
exist for \emph{any} $\tau \in R$ and any $\phi \in C^1_c(\Omega)$, $\bfphi \in C^1_c(\Omega; R^3)$. In addition, 
\begin{align}
\left< \vr(\tau-, \cdot); \phi \right> &= \left< \vr(\tau+, \cdot); \phi \right>, \ \mbox{and}\ 
\tau \mapsto \left< \vr(\tau, \cdot); \phi \right> \in BC(R), \br
\left< \vm(\tau-, \cdot) ; \bfphi \right> &= \left< \vm(\tau+, \cdot) ; \bfphi \right>, \ \mbox{and}\ 
\tau \mapsto \left< \vm(\tau, \cdot) ; \bfphi \right> \in BC(R), 
\nonumber
\end{align}
and 
\begin{align} 
\left< S(\tau-, \cdot); \phi \right> &\leq \left< S(\tau+, \cdot); \phi \right> \ \mbox{whenever}\ \phi \geq 0, \br
\mbox{and}\
\tau &\mapsto \left< S(\tau-, \cdot); \phi \right> = h_\phi (\tau) + g_\phi (\tau)	,\ g_\phi \in C_{\rm loc}(R),\ 
h_\phi \ \mbox{non--decreasing c\` agl\` ad.}
	\label{L1}
	\end{align}

The trajectory space can be therefore identified with ``weakly c\` agl\` ad'' and bounded functions defined on $R$. To this end, consider the Hilbert space
\[
W^{k,2}_0(\Omega),\ k > \frac{3}{2} \ \mbox{so that}\ W^{k,2}_0 \hookrightarrow C(\Ov{\Omega})
\]
with an orthonormal basis of smooth functions $\{ \phi_n \}_{n=1}^\infty$. Similarly, we consider the same space of vector valued 
functions $W^{1,2}_0(\Omega; R^3)$ with a basis $\{ \bfphi_n \}_{n=1}^\infty$.
Finally, we define a metrics 
\begin{align} 
d_{\mathcal{T}} &\left[ (\vr^1, S^1, \vm^1) ; (\vr^2, S^2, \vm^2) \right] \br 
&= \sum_{n=1}^\infty \frac{1}{2^n} \int_{- \infty}^{\infty} \exp \left(-t^2 \right) G \left( \| \left< \vr^1 - \vr^2; \phi_n \right> \|_{C[-t,t]} \right) \dt \br &+ 
\sum_{n=1}^\infty \frac{1}{2^n} \int_{- \infty}^{\infty} \exp \left(-t^2 \right) G \left( \| \left< \vm^1 - \vm^2; \bfphi_n \right> \|_{C([-t,t]; R^3)} \right) \dt \br 
&= \sum_{n=1}^\infty \frac{1}{2^n} \int_{- \infty}^{\infty} \exp \left(-t^2 \right) G \left( \left[ \left< S^1; \phi_n \right>; 
\left< S^2; \phi_n \right> \right]_{D[-t,t]} \right) \dt, 
 \label{L2}	
	\end{align}
where 
\[
G(Z) = \frac{Z}{1 + Z} 
\]
and $D[-t,t]$ denotes the Skorokhod space of of c\` agl\` ad functions defined on $[-t,t]$ with the associated 
complete metrics $[ \cdot ; \cdot ]_{D[-t,t]}$, see e.g. Whitt \cite{Whitt}.

The trajectory space is defined as 
\[
\mathcal{T} = \cup_{L = 1}^\infty \mathcal{T}_L, 
\]
where 
\begin{align}
\mathcal{T}_L = \Big\{ (\vr, S, \vm) \ \Big| \ &\vr \in L^\infty (R; W^{-k,2}(\Omega)), \ \left< \vr; \phi_n \right> \in C(R),\ n=1,2,\dots, \br  &\sup_{t \in R} \| \vr(t, \cdot) \|_{W^{-k,2}(\Omega)} \leq L, \br 
&\vm \in L^\infty (R; W^{-k,2}(\Omega; R^3)), \ \left< \vm; \bfphi_n  \right> \in C(R),\ n=1,2,\dots, \br
&\sup_{t \in R} \| \vm(t, \cdot) \|_{W^{-k,2}(\Omega; R^3)} \leq L, \br 
&S \in L^\infty (R; W^{-k,2}(\Omega)),\ \left< S; \phi_n \right> \ \mbox{c\` agl\` ad in}\ R, \ n=1,2,\dots, \br
&\sup_{t \in R} \| S(t, \cdot) \|_{W^{-k,2}(\Omega)} \leq L \Big\}. 
\nonumber
\end{align}
Note that the trajectory space is larger then the set of entire solutions to the Navier--Stokes--Fourier system 
and consists of time dependent functionals ranging in the space of distributions on $\Omega$.

Each set $\mathcal{T}_L$ endowed with the metrics $d_{\mathcal{T}}$ is a Polish space. We define inductive topology on~$\mathcal{T}$:
\begin{align}
(\vr_n, S_n, \vm_n) \to (\vr, S, \vm) \ \ \mbox{in} \ \mathcal{T} \ \Leftrightarrow \ & \mbox{{\bf (i)} there exists}\ L \ \mbox{such that}\ (\vr_n, S_n, \vm_n) \in \mathcal{T}_L
\ \mbox{for all}\ n=1,2,\dots \br 
&\mbox{{\bf (ii)}}\  d_{\mathcal{T}} \left[ (\vr_n, S_n, \vm_n); (\vr, S, \vm) \right] \to 0. 
\nonumber
\end{align}
Our choice of the topology of the trajectory space may seem a bit awkward at the first glance but accommodates 
the instantaneous convergence of the state variables. Alternatively, a weaker $L^p$ topology can be used being equivalent on the attractor $\mathcal{A}$ to $d_{\mathcal{T}}$.

\subsection{Attractor}

As the weak solutions are not (known to be) uniquely determined by the initial/boundary data, we adopt the approach of Sell \cite{SEL} and M\' alek and Ne\v cas \cite{MANE} 
replacing the standard phase space by the trajectory space $\mathcal{T}$. Here and hereafter, we always assume that the principal hypotheses of Theorem \ref{MT1} concerning the constitutive 
relations are satisfied. Moreover, we fix the total mass of the fluid, 
\begin{equation} \label{L3}
	\intO{ \vr(t, \cdot) } = M > 0.
\end{equation}
Accordingly, the set $\mathcal{A}$, 
\begin{align} 
	\mathcal{A} = \Big\{ &(\vr, S, \vm) \ \Big| \ (\vr, S, \vm) \ \mbox{a weak solution of the Navier--Stokes--Fourier system} \br &\mbox{in the sense of Definition \ref{Dw1} 
	on the time interval}\ t \in R \Big\}, 	
	\label{L4}
	\end{align}
is a natural candidate to be \emph{global attractor} in the trajectory space $\mathcal{T}$. Indeed, as shown in Theorem~\ref{MT1}, 
\begin{equation} \label{L5}
\sup_{t \in R} \intO{ E \left( \vr, S, \vm \right) (t,\cdot) } \leq \mathcal{E}_{\infty} < \infty.
\end{equation}
In particular, $\mathcal{A} \subset \mathcal{T}_L \subset \mathcal{T}$ for a sufficiently large $L$.

\begin{Lemma} \label{LL1}
	Under the hypotheses of Theorem \ref{MT1}, the set $\mathcal{A}$ is 
	\begin{itemize} 
		\item non--empty; 
		\item time--shift invariant, 
		\[
		(\vr, S, \vm) \in \mathcal{A}\ \Rightarrow \ (\vr, S, \vm) (\cdot + T) \in \mathcal{A} \ \mbox{for any}\ T \in R;
		\]
	\item compact in the metric topology $(\mathcal{T}_L, d_{\mathcal{T}})$ for a sufficiently large $L$.
	
	\end{itemize}
	
	\end{Lemma}

\begin{proof}
	
	As shown in \cite[Theorem 4.2]{ChauFei}, the Navier--Stokes--Fourier system with the boundary conditions \eqref{i5}, \eqref{i6} admits a global--in--time weak solution $(\vr, \vt, \vu)$ on the time interval 
	$[0, \infty)$ for any initial data with finite energy. It follows from Theorem \ref{MT2} that there exists a sequence of times 
	$T_n \to \infty$ such that $((\vr, \vt, \vu)(\cdot + T_n))_{n = 1}^\infty$ converge to a weak solution $(\tvr, \tvt, \tvu)$ of the same problem defined for all $t \in R$ with globally bounded energy. Obviously, 
	\[
	(\vr, S, \vm) = (\tvr,  \tvr s(\tvr, \tvt), \tvr \tvu) \in \mathcal{A}.
	\]
	Moreover, as the underlying system is autonomous, the set $\mathcal{A}$ is time--shift  invariant.
	
	Compactness of the set $\mathcal{A}$ follows again from Theorem \ref{MT2}, where we consider $T_n = - \infty$. As pointed out, $\mathcal{A} \subset \mathcal{T}_L$, where the latter is a metric space; whence compactness is equivalent to sequential compactness. At the level of the density, convergence in the metric $d_{\mathcal{T}}$ follows  from \eqref{M1}. Moreover, since the 
	momenta $\vr \vu$ satisfy equation \eqref{w5}, they are precompact in the topology of the space 
	\[
	C_{\rm weak, loc} (R; L^{\frac{5}{4}}(\Omega; R^3)), 
	\]
	which implies compactness in the momentum component of $\mathcal{T}_L$ with the $d_{\mathcal{T}}$ metrics. 
	
	Thus it remains to show compactness at the level of the total entropy $S$. First observe that 
	\[
	\tau \mapsto \intO{ S (\tau, \cdot) \phi } - \int_0^\tau \intO{ \vr s(\vr, \vt) \vu \cdot \Grad \phi } \dt + \int_0^\tau \intO{ \frac{\kappa (\vt) }{\vt} \cdot \Grad \phi } \dt 
\]	
is a non--decreasing function of $\tau$ for any test function $\phi \in C^1_c(\Omega)$, $\phi \geq 0$. In view of boundedness of the total energy, the family 
\[
\tau \mapsto \int_0^\tau \intO{ \vr s(\vr, \vt) \vu \cdot \Grad \phi } \dt + \int_0^\tau \intO{ \frac{\kappa (\vt) }{\vt} \cdot \Grad \phi } \dt 
\] 
is precompact in $C_{\rm loc}(R)$; whence precompactness of $S$ in $d_{\mathcal{T}}$ reduces to precompactness of a \emph{non--decreasing (in time)} sequence of functions
\[
\tau \mapsto \left< \tS_n ; \phi \right> \equiv \intO{ S_n (\tau, \cdot) \phi } - \int_0^\tau \intO{ \vr_n s(\vr_n, \vt_n) \vu_n \cdot \Grad \phi } \dt + \int_0^\tau \intO{ \frac{\kappa (\vt_n) }{\vt_n} \cdot \Grad \phi } \dt
\]
with respect to the metrics $d$ of the Skorokhod space $D[-M,M]$ of c\` agl\` ad functions defined on compact time intervals $[-M, M]$. To this end, we recall the criterion due to Whitt \cite[Chapter 12, Corollary 12.5.1]{Whitt}:
	\medskip
	
	\hrule
	
	\medskip
		                                   
		Let 
		\[
		h_n:[-M,M] \to R
        \]                                               
		be a sequence of monotone functions. 
		
		\medskip
		
		Then 
		\[
		d[h_n; h] \to 0 \ \mbox{as}\ n \to \infty \ \mbox{for some}\ h \in D[-M,M] 
		\]
		if and only if 
		\[
		h_n(t) \to h(t) \ \mbox{for all} \ t \ \mbox{belonging to a dense set in}\ [-M,M]\ \mbox{including the end points}\  -M \ \mbox{and}\  M.
		\]
		
\medskip \hrule \medskip

Now, as $S_n$ are the total entropies generated by a family of uniformly bounded weak solutions of the Navier--Stokes--Fourier system, we have 
\begin{align} 
\left< \tS_n (\tau, \cdot); \phi \right> &\to \left< \tS (\tau, \cdot); \phi \right> \br
& \equiv
\intO{ S (\tau, \cdot) \phi } - \int_0^\tau \intO{ \vr s(\vr, \vt) \vu \cdot \Grad \phi } \dt + \int_0^\tau \intO{ \frac{\kappa (\vt) }{\vt} \cdot \Grad \phi } \dt
\nonumber
\end{align}
for a.a. $\tau \in R$ at least for a suitable subsequence, where $(\vr, S = \vr s(\vr, \vt), \vm = \vr \vu) \in \mathcal{A}$, is another entire solution of the same problem.
In particular, 
\[
\left< \tS_n (\tau, \cdot); \phi \right> \to \left< \tS (\tau, \cdot); \phi \right> \ \mbox{for a dense set of times for any compact interval} \ [-M, M];
\]
whence, in accordance with the above convergence criterion,                                                       
\[                                                                                                                                                                   
\left< \tS_n (\tau, \cdot); \phi \right> \to \left< \tS (\tau, \cdot); \phi \right> \ \mbox{in} \ D[-M, M] \ \mbox{for a.a.}\ M > 0                                                                                                                      
\]
yielding the desired conclusion
\[                                                                                                                                                                   
\left< S_n (\tau, \cdot); \phi \right> \to \left< S (\tau, \cdot); \phi \right> \ \mbox{in} \ D[-M, M] \ \mbox{for a.a.}\ M > 0.                                                                                                                      
\]

\end{proof}
 
In accordance with Theorems \ref{MT1}, \ref{MT2}, and Lemma \ref{LL1}, we may state our main result concerning the existence of a trajectory attractor for the Rayleigh--B\' enard problem:

\begin{mdframed}[style=MyFrame]
	
	\begin{Theorem} [{\bf Trajectory attractor}] \label{LT1}
	
	Let $M > 0$, $\mathcal{E}_0$ be given. Let $\mathcal{F}[M, \mathcal{E}_0]$ be a family of weak solutions to the Rayleigh--B\' enard problem for the Navier--Stokes--Fourier system
	on the time interval $(0, \infty)$
	 satisfying 
	\[
	\intO{ \vr } = M,\ {\rm ess} \limsup_{\tau \to 0+} \intO{ E(\vr, S, \vm) (\tau, \cdot) } \leq \mathcal{E}_0.
	\]
	We identify the set $\mathcal{F}[M, \mathcal{E}_0]$ with a subset of the trajectory space $\mathcal{T}$ extending 
	\[ 
	\vr(\tau, \cdot) = \lim_{t \to 0+} \vr(t, \cdot),\ \vm (\tau, \cdot) = \lim_{t \to 0+} \vm(t, \cdot),\ 
	0 \leq S(\tau, \cdot) \leq \lim_{t \to 0+} S(t, \cdot) \ \mbox{for}\ \tau < 0,
	\]
	where the limits are understood in the weak (distributional) sense.
	
	Then for any $\ep > 0$, there exists a time $T(\ep)$ such that 
	\[
	d_{\mathcal{T}} [ (\vr, S, \vm)(\cdot + T); \mathcal{A} ] < \ep \ \mbox{for any}\ (\vr, S, \vm) \in \mathcal{F}[M, \mathcal{E}_0] 
	\ \mbox{and any}\ T > T(\ep).
	\]
	
	\end{Theorem}
	
	\end{mdframed}

\subsection{Stationary statistical solutions}

Following the ideas of the preceding section we are ready to identify {\it statistical  solutions} with shift invariant probability measures on  ${\mathcal T}$, which are supported by solutions to the Navier-Stokes-Fourier system. In accordance with Theorem \ref{LT1}, these shift invariant probability measures supported by the global trajectory attractor $\mathcal{A}$.

The construction of a bounded invariant measure is the same as in \cite{FanFeiHof}, note that a similar approach in the incompressible setting was used by Foias, Rosa and Temam \cite{FoRoTe2}, \cite{FoRoTe1}.
Given a trajectory $(\vr, S, \vm) \in \mathcal{T}$ we consider a probability measure 
\[
\mathcal{V}_T \equiv \frac{1}{T} \int_0^T \delta_{(\vr, S, \vm) (\cdot + t)} \dt, 
\]
where $\delta$ denotes the Dirac mass. Obviously, $\mathcal{V}_T$ is a probability measure. If, in addition, 
\[
(\vr, S, \vm) \in \mathcal{A},
\] 
then $\mathcal{V}_T \in \mathfrak{P} (\mathcal{A})$, where $\mathfrak{P}(\mathcal{A})$ denotes the set of all probablity measures on a compact Polish space $\mathcal{A}$. In particular, the family
\[
\{ \mathcal{V}_T \}_{T \geq 0} \ \mbox{is tight.}
\]
By Prokhorov theorem, there is a sequence $T_n \to \infty$ such that 
\[
\mathcal{V}_{T_n} \to \mathcal{V} \ \mbox{narowly in}\ \mathfrak{P}(\mathcal{A}). 
\]
Finally, exactly as in \cite[Section 5.1]{FanFeiHof}, we may show that the measure $\mathcal{V}$ is time--shift  invariant, meaning 
\begin{equation} \label{L6}
\mathcal{V}[ \mathfrak{B}(\cdot + T) ] = \mathcal{V}[ \mathfrak{B} ] \ \mbox{for any Borel set}\ \mathfrak{B} \subset \mathcal{T}.
\end{equation}
A Borel probability measure $\mathcal{V} \in \mathfrak{P}(\mathcal{A})$ enjoying the property \eqref{L6} is called \emph{statistical stationary solution} of the 
Rayleigh--B\' enard problem for the Navier--Stokes--Fourier system.

Finally, observe that the above construction may be restricted to any shift--invariant subset $\mathcal{U} \subset \mathcal{A}$.

\begin{mdframed}[style=MyFrame]

\begin{Theorem}
Let ${\mathcal U}\subset \mathcal{A}$ be a non-empty time--shift  invariant set, meaning 
\[
(\vr, S, \vm) \in \mathcal{U} \ \Rightarrow \ (\vr, S, \vm)(\cdot + T) \in \mathcal{U} \ \mbox{for any}\ T \in R.
\]

Then there exists a stationary statistical solution ${\mathcal V}$ supported by $\Ov{\mathcal{U}}$:
\begin{itemize}
\item ${\mathcal V}$ is a Borel probability measure, $\mathcal{V} \in \mathfrak{P}(\Ov{\mathcal{U}})$;
\item  $\supp {\mathcal V}\subset \overline{\mathcal U}$, where the closure of a ${\mathcal U}$ is a compact invariant set;
\item ${\mathcal V}$ is shift invariant, i.e., ${\mathcal V}[\mathfrak{B}]={\mathcal V}[{\mathfrak{B}(\cdot+T)}]$
for any Borel set $\mathfrak{B}\subset {\mathcal T}$ and any $T\in R$.
\end{itemize}

\end{Theorem}

\end{mdframed}

\subsection{Convergence of ergodic means}

We conclude this section by a direct application of Birkhoff--Khinchin ergodic theorem. Similarly to \cite[Section 5]{FanFeiHof}, we may consider the state space 
\[
H = W^{-k,2}(\Omega) \times W^{-k,2}(\Omega) \times W^{-k,2}(\Omega; R^3).
\]
Let $\mathcal{V} \in \mathfrak{P}( \mathcal{T})$ be a statistical stationary solution, and thus a Borel probability measure. Consider   a probability basis $(\mathcal{T}, \mathfrak{B}[\mathcal{T}], \mathcal{V} )$ where $\mathcal{T}$ is the trajectory space and $\mathfrak{B}[ \mathcal{T} ]$ is the family of Borel sets.
%Given a trajectory $(\vr, S, \vm) \in \mathcal{T}$, we may consider the associated canonical process 
%\[
%(\vr, S, \vm) \in \mathcal{T}  \times \tau \in  R \mapsto (\vr, S, \vm) (\tau, \cdot) \in H. 
%\]
Given a trajectory $(\vr, S, \vm) \in \mathcal{T}$ and 
$\tau\in R$ we may consider the associated canonical process 
\[
(\vr, S, \vm)   \times \tau  \mapsto (\vr, S, \vm) (\tau, \cdot) \in H. 
\]

As $\mathcal{V}$ is shift invariant, the above process is a stationary process with respect to the probability basis 
$(\mathcal{T}, \mathfrak{B}[\mathcal{T}], \mathcal{V} )$ defined for $\tau \in R$, with c\` agl\` ad paths ranging in $H$.

Exactly as in \cite[Theorem 6.4, Section 6]{FanFeiHof} we can establish the following result.

\begin{mdframed}[style=MyFrame]
	
	\begin{Theorem} [{\bf Convergence of ergodic averages}] \label{LT2}
		
		Let $\mathcal{V}$ be a stationary statistical solution and $(\vr, S, \vm)$ the associated stationary process. Let $F: H \to R$ be a Borel measurable function such that 
		\[
		\int_{\mathcal{T}} |F (\vr(0, \cdot), S(0, \cdot), \vm (0, \cdot) | \ \D \mathcal{V} < \infty.
		\]
		
		Then there exists a measurable function $\Ov{F}$, 
		\[
		\Ov{F}: (\mathcal{T}, \mathcal{V} ) \to R 
		\]
		such that
		\[
		\frac{1}{T} \int_0^T F(\vr(t, \cdot), S(t, \cdot), \vm(t, \cdot) ) \dt \to \Ov{F} \ \mbox{as}\ T \to \infty
		\]
$\mathcal{V}-$a.s. and in $L^1 (\mathcal{T}, \mathcal{V} )$.

		\end{Theorem}
	
	\end{mdframed}

\def\cprime{$'$} \def\ocirc#1{\ifmmode\setbox0=\hbox{$#1$}\dimen0=\ht0
	\advance\dimen0 by1pt\rlap{\hbox to\wd0{\hss\raise\dimen0
			\hbox{\hskip.2em$\scriptscriptstyle\circ$}\hss}}#1\else {\accent"17 #1}\fi}

%\bibliography{citace}

\begin{thebibliography}{10}
	
	\bibitem{BEROFO}
	S.~E. Bechtel, F.J. Rooney, and M.G. Forest.
	\newblock Connection between stability, convexity of internal energy, and the
	second law for compressible {N}ewtonian fuids.
	\newblock {\em J. Appl. Mech.}, {\bf 72}:299--300, 2005.
	
	\bibitem{BEL1}
	F.~Belgiorno.
	\newblock Notes on the third law of thermodynamics, {I}.
	\newblock {\em J. Phys. A}, {\bf 36}:8165--8193, 2003.
	
	\bibitem{BEL2}
	F.~Belgiorno.
	\newblock Notes on the third law of thermodynamics, ii.
	\newblock {\em J. Phys. A}, {\bf 36}:8195--8221, 2003.
	
	\bibitem{Borm2}
	A.~S. Bormann.
	\newblock The onset of convection in the {R}ayleigh-{B}enard problem for
	compressible fluids.
	\newblock {\em Continuum Mech. Thermodyn.}, {\bf 13}:9--23, 2001.
	
	\bibitem{Borm1}
	A.~S. Bormann.
	\newblock Numerical linear stability analysis for compressible fluids.
	\newblock In {\em Analysis and numerics for conservation laws}, pages 93--105.
	Springer, Berlin, 2005.
	
	\bibitem{BRDE}
	D.~Bresch and B.~Desjardins.
	\newblock Stabilit{\' e} de solutions faibles globales pour les {\' e}quations
	de {N}avier-{S}tokes compressibles avec temp{\' e}rature.
	\newblock {\em C.R. Acad. Sci. Paris}, {\bf 343}:219--224, 2006.
	
	\bibitem{BreJab}
	D.~Bresch and P.-E. Jabin.
	\newblock Global existence of weak solutions for compressible {N}avier-{S}tokes
	equations: thermodynamically unstable pressure and anisotropic viscous stress
	tensor.
	\newblock {\em Ann. of Math. (2)}, {\bf 188}(2):577--684, 2018.
	
	\bibitem{CaJoTuWh}
	Y.~Cao, M.~S. Jolly, E.~S. Titi, and J.~P. Whitehead.
	\newblock Algebraic bounds on the {R}ayleigh-{B}\'{e}nard attractor.
	\newblock {\em Nonlinearity}, {\bf 34}(1):509--531, 2021.
	
	\bibitem{ChauFei}
	N.~Chaudhuri and E.~Feireisl.
	\newblock {N}avier--{S}tokes--{F}ourier system with {D}irichlet boundary
	conditions.
	\newblock {\em {\bf arxiv preprint No. 2106.05315}}, 2021.
	
	\bibitem{ChoNobOtt}
	A.~Choffrut, C.~Nobili, and F.~Otto.
	\newblock Upper bounds on {N}usselt number at finite {P}randtl number.
	\newblock {\em J. Differential Equations}, {\bf 260}(4):3860--3880, 2016.
	
	\bibitem{CoFoMaTe}
	P.~Constantin, C.~Foias, O.~P. Manley, and R.~Temam.
	\newblock Determining modes and fractal dimension of turbulent flows.
	\newblock {\em J. Fluid Mech.}, {\bf 150}:427--440, 1985.
	
	\bibitem{DBPB}
	K.~E. Daniels, O.~Brausch, W.~Pesch, and E.~Bodenschatz.
	\newblock Competition and bistability of ordered undulations and undulation
	chaos in inclined layer convection.
	\newblock {\em J. Fluid Mech.}, {\bf 597}:261--282, 2008.
	
	\bibitem{DAVI}
	P.~A. Davidson.
	\newblock {\em Turbulence:{A}n introduction for scientists and engineers}.
	\newblock Oxford University Press, Oxford, 2004.
	
	\bibitem{FanFeiHof}
	F.~Fanelli, E.~Feireisl, and M.~Hofmanov\'{a}.
	\newblock Ergodic theory for energetically open compressible fluid flows.
	\newblock {\em Phys. D}, {\bf 423}:Paper No. 132914, 25, 2021.
	
	\bibitem{FeiKwo}
	E.~Feireisl and Y.-S. Kwon.
	\newblock Asymptotic stability of solutions to the {N}avier-{S}tokes-{F}ourier
	system driven by inhomogeneous {D}irichlet boundary conditions.
	\newblock {\em Archive Preprint Series}, 2021.
	\newblock {\bf arxiv preprint No. 2109.00980}.
	
	\bibitem{FeNo6A}
	E.~Feireisl and A.~Novotn\'y.
	\newblock {\em Singular limits in thermodynamics of viscous fluids}.
	\newblock Advances in Mathematical Fluid Mechanics. Birkh\"auser/Springer,
	Cham, 2017.
	\newblock Second edition.
	
	\bibitem{FeiNov20}
	E.~Feireisl and A.~Novotn\'{y}.
	\newblock Navier-{S}tokes-{F}ourier {S}ystem with {G}eneral {B}oundary
	{C}onditions.
	\newblock {\em Comm. Math. Phys.}, {\bf 386}(2):975--1010, 2021.
	
	\bibitem{FP9}
	E.~Feireisl and H.~Petzeltov{\'a}.
	\newblock Large-time behaviour of solutions to the {N}avier-{S}tokes equations
	of compressible flow.
	\newblock {\em Arch. Rational Mech. Anal.}, {\bf 150}:77--96, 1999.
	
	\bibitem{FP14}
	E.~Feireisl and H.~Petzeltov{\'a}.
	\newblock Bounded absorbing sets for the {N}avier-{S}tokes equations of
	compressible fluid.
	\newblock {\em Commun. Partial Differential Equations}, {\bf 26}:1133--1144,
	2001.
	
	\bibitem{FP20}
	E.~Feireisl and H.~Petzeltov\'{a}.
	\newblock On the long-time behaviour of solutions to the
	{N}avier-{S}tokes-{F}ourier system with a time-dependent driving force.
	\newblock {\em J. Dynam. Differential Equations}, {\bf 19}(3):685--707, 2007.
	
	\bibitem{FeiPr}
	E.~Feireisl and D.~Pra{\v z}{\' a}k.
	\newblock {\em Asymptotic behavior of dynamical systems in fluid mechanics}.
	\newblock AIMS, Springfield, 2010.
	
	\bibitem{FoMaTe}
	C.~Foias, O.~Manley, and R.~Temam.
	\newblock Attractors for the {B}\'{e}nard problem: existence and physical
	bounds on their fractal dimension.
	\newblock {\em Nonlinear Anal.}, {\bf 11}(8):939--967, 1987.
	
	\bibitem{FoRoTe2}
	C.~Foias, R.~M.~S. Rosa, and R.~M. Temam.
	\newblock Convergence of time averages of weak solutions of the
	three-dimensional {N}avier-{S}tokes equations.
	\newblock {\em J. Stat. Phys.}, {\bf 160}(3):519--531, 2015.
	
	\bibitem{FoRoTe1}
	C.~Foias, R.~M.~S. Rosa, and R.~M. Temam.
	\newblock Properties of stationary statistical solutions of the
	three-dimensional {N}avier-{S}tokes equations.
	\newblock {\em J. Dynam. Differential Equations}, {\bf 31}(3):1689--1741, 2019.
	
	\bibitem{GALN}
	G.~P. Galdi.
	\newblock {\em An introduction to the mathematical theory of the {N}avier -
		{S}tokes equations, Second Edition}.
	\newblock Springer-Verlag, New York, 2003.
	
	\bibitem{GEHEHI}
	M.~Gei{\ss}ert, H.~Heck, and M.~Hieber.
	\newblock On the equation {${\rm div}\,u=g$} and {B}ogovski\u\i's operator in
	{S}obolev spaces of negative order.
	\newblock In {\em Partial differential equations and functional analysis},
	volume 168 of {\em Oper. Theory Adv. Appl.}, pages 113--121. Birkh\"auser,
	Basel, 2006.
	
	\bibitem{MANE}
	J.~M{\' a}lek and J.~Ne{\v c}as.
	\newblock A finite-dimensional attractor for the three dimensional flow of
	incompressible fluid.
	\newblock {\em J. Differential. Equations}, {\bf 127}:498--518, 1996.
	
	\bibitem{MANI1}
	A.~Matsumura and T.~Nishida.
	\newblock The initial value problem for the equations of motion of viscous and
	heat-conductive gases.
	\newblock {\em J. Math. Kyoto Univ.}, {\bf 20}:67--104, 1980.
	
	\bibitem{MANI}
	A.~Matsumura and T.~Nishida.
	\newblock The initial value problem for the equations of motion of compressible
	and heat conductive fluids.
	\newblock {\em Comm. Math. Phys.}, {\bf 89}:445--464, 1983.
	
	\bibitem{NobOtt1}
	C.~Nobili and F.~Otto.
	\newblock Limitations of the background field method applied to
	{R}ayleigh-{B}\'{e}nard convection.
	\newblock {\em J. Math. Phys.}, {\bf 58}(9):093102, 46, 2017.
	
	\bibitem{NovPok07}
	A.~Novotn\'{y} and M.~Pokorn\'{y}.
	\newblock Stabilization to equilibria of compressible {N}avier-{S}tokes
	equations with infinite mass.
	\newblock {\em Comput. Math. Appl.}, 53(3-4):437--451, 2007.
	
	\bibitem{NOS1}
	A.~Novotn{\' y} and I.~Stra{\v s}kraba.
	\newblock Stabilization of weak solutions to compressible {N}avier-{S}tokes
	equations.
	\newblock {\em J. Math. Kyoto Univ.}, {\bf 40}:217--245, 2000.
	
	\bibitem{NOST}
	A.~Novotn{\' y} and I.~Stra{\v s}kraba.
	\newblock Convergence to equilibria for compressible {N}avier-{S}tokes
	equations with large data.
	\newblock {\em Annali Mat. Pura Appl.}, {\bf 169}:263--287, 2001.
	
	\bibitem{OTTOetal1}
	F.~Otto, S.~Pottel, and C.~Nobili.
	\newblock Rigorous bounds on scaling laws in fluid dynamics.
	\newblock In {\em Mathematical thermodynamics of complex fluids}, volume 2200
	of {\em Lecture Notes in Math.}, pages 101--145. Springer, Cham, 2017.
	
	\bibitem{SEL}
	G.~R. Sell.
	\newblock Global attractors for the three-dimensional {N}avier-{S}tokes
	equations.
	\newblock {\em J. Dynamics Differential Equations}, {\bf 8}(1):1--33, 1996.
	
	\bibitem{VAZA}
	A.~Valli and M.~Zajaczkowski.
	\newblock {N}avier-{S}tokes equations for compressible fluids: Global existence
	and qualitative properties of the solutions in the general case.
	\newblock {\em Commun. Math. Phys.}, {\bf 103}:259--296, 1986.
	
	\bibitem{WanWan}
	X.~Wang and W.~Wang.
	\newblock On global behavior of weak solutions to the {N}avier-{S}tokes
	equations of compressible fluid for {$\gamma=5/3$}.
	\newblock {\em Bound. Value Probl.}, pages 2015:176, 13, 2015.
	
	\bibitem{Whitt}
	W.~Whitt.
	\newblock {\em Stochastic-process limits}.
	\newblock Springer Series in Operations Research. Springer-Verlag, New York,
	2002.
	\newblock An introduction to stochastic-process limits and their application to
	queues.
	
\end{thebibliography}

%\bibliographystyle{plain}

\end{document}